\documentclass{article}

\usepackage{amssymb,amsmath,amsthm,mathrsfs}
\usepackage{verbatim}
\usepackage[a4paper,includeheadfoot,left=30mm,right=30mm,top=25mm,bottom=30mm]{geometry} 

\newtheorem{thm}{Theorem}[section]
\newtheorem{lem}[thm]{Lemma}

\newtheorem{defn}[thm]{Definition}

\theoremstyle{definition}
\newtheorem{cond}[thm]{Condition}

\newtheorem{ntn}[thm]{Notation}
\theoremstyle{remark}
\newtheorem{rmk}[thm]{Remark}
\newtheorem{eg}[thm]{Example}

\newcommand{\TextNorm}[1]{\textrm{\textmd{\textup{#1}}}}
\newcommand{\hsforall}{\hspace{1mm}\forall\hspace{1mm}}						 
\newcommand{\hsexists}{\hspace{1mm}\exists\hspace{1mm}} 					 
\newcommand{\res}{\operatorname*{Res}}                             
\renewcommand{\Re}{\operatorname*{Re}}                             
\renewcommand{\Im}{\operatorname*{Im}}                             
\renewcommand{\d}{\ensuremath{\,\mathrm{d}}}							         
\newcommand{\T}{\text{\TextNorm{T}}\hspace{0.1mm}}							   
\newcommand{\Ymax}{Y_{\mathrm{max}}}                               
\newcommand{\DeltaP}{\Delta_{\mathrm{PDE}}\hspace{0.5mm}}          
\newcommand{\DeltaPprime}{\Delta'_{\mathrm{PDE}}\hspace{0.5mm}}    
\newcommand{\M}[3]{#1_{#2\hspace{0.5mm}#3}}                        
\renewcommand{\geq}{\geqslant}                                     
\renewcommand{\leq}{\leqslant}                                     
\newcommand{\BE}{\begin{equation}}                                 
\newcommand{\EE}{\end{equation}}                                   
\newcommand{\BES}{\begin{equation*}}                               
\newcommand{\EES}{\end{equation*}}                                 
\newcommand{\BP}{\begin{pmatrix}}                                  
\newcommand{\EP}{\end{pmatrix}}                                    

\newcommand{\superscript}[1]{\ensuremath{^{\textrm{#1}}}}

\newcommand{\Thns}[0]{\superscript{th}}
\newcommand{\Th}[0]{\Thns~}

\def\clap#1{\hbox to 0pt{\hss#1\hss}}
\def\mathclap{\mathpalette\mathclapinternal}
\def\mathclapinternal#1#2{\clap{$\mathsurround=0pt#1{#2}$}}

\numberwithin{equation}{section}

\title{Well-posed two-point initial-boundary value \\ problems with arbitrary boundary conditions}
\author{David A. Smith \\ Department of Mathematics, University of Reading RG\textup{6} \textup{6}AX \\ email\textup{: \texttt{d.a.smith@reading.ac.uk}}}

\begin{document}
\maketitle

\abstract{We study initial-boundary value problems for linear evolution equations of arbitrary spatial order, subject to arbitrary linear boundary conditions and posed on a rectangular 1-space, 1-time domain. We give a new characterisation of the boundary conditions that specify well-posed problems using Fokas' transform method. We also give a sufficient condition guaranteeing that the solution can be represented using a series.

The relevant condition, the analyticity at infinity of certain meromorphic functions within particular sectors, is significantly more concrete and easier to test than the previous criterion, based on the existence of admissible functions.}

\section{Introduction} \label{sec:P1:Intro}

In this work, we consider
\smallskip

\noindent{\bfseries The initial-boundary value problem $\Pi(n,A,a,h,q_0)$:}
Find $q\in C^\infty([0,1]\times[0,T])$ which satisfies the linear, evolution, constant-coefficient partial differential equation
\BE \label{eqn:P1:Intro:PDE}
\partial_tq(x,t) + a(-i\partial_x)^nq(x,t) = 0
\EE
subject to the initial condition
\BE \label{eqn:P1:Intro:IC}
q(x,0) = q_0(x)
\EE
and the boundary conditions
\BE \label{eqn:P1:Intro:BC}
A\left(\partial_x^{n-1}q(0,t),\partial_x^{n-1}q(1,t),\partial_x^{n-2}q(0,t),\partial_x^{n-2}q(1,t),\dots,q(0,t),q(1,t)\right)^\T = h(t),
\EE

\noindent where the pentuple $(n,A,a,h,q_0)\in\mathbb{N}\times\mathbb{R}^{n\times2n}\times\mathbb{C}\times(C^\infty[0,T])^n \times C^\infty[0,1]$ is such that
\begin{description}
\item[$(\Pi1)$]{the \emph{order} $n\geq 2$,}
\item[$(\Pi2)$]{the \emph{boundary coefficient matrix} $A$ is in reduced row-echelon form,}
\item[$(\Pi3)$]{if $n$ is odd then the \emph{direction coefficient} $a=\pm i$, if $n$ is even then $a=e^{i\theta}$ for some $\theta\in[-\pi/2,\pi/2]$,}
\item[$(\Pi4)$]{the \emph{boundary data} $h$ and the \emph{initial datum} $q_0$ are compatible in the sense that
\BE \label{eqn:P1:Intro:Compatibility}
A\left(q_0^{(n-1)}(0),q_0^{(n-1)}(1),q_0^{(n-2)}(0),q_0^{(n-2)}(1),\dots,q_0(0),q_0(1)\right)^\T = h(0).
\EE}
\end{description}

Provided $\Pi$ is well-posed, in the sense of admitting a unique, smooth solution, its solution may be found using Fokas' unified transform method~\cite{Fok2008a,FP2001a}. The representation thus obtained is a contour integral of transforms of the initial and boundary data. Certain problems, for example those with periodic boundary conditions, may be solved using classical methods such as Fourier's separation of variables~\cite{Fou1822a}, to yield a representation of the solution as a discrete Fourier series. By the well-posedness of $\Pi$, these are two different representations of the same solution.

For individual examples, Pelloni~\cite{Pel2005a} and Chilton~\cite{Chi2006a} discuss a method of recovering a series representation from the integral representation through a contour deformation and a residue calculation. Particular examples have been identified of well-posed problems for which this deformation fails but there is no systematic method of determining its applicability.

Pelloni~\cite{Pel2004a} uses Fokas' method to decide well-posedness of a class of problems with uncoupled, non-Robin boundary conditions giving an explicit condition, the number that must be specified at each end of the space interval, whose validity may be ascertained immediately. However there exist no criteria for well-posedness that are at once more general than Pelloni's and simpler to check than the technical `admissible set' characterisation of~\cite{FP2001a}.

The principal result of this work is a new characterisation of well-posedness. The condition is the decay of particular integrands within certain sectors of the complex plane. Indeed, let $D=\{\rho\in\mathbb{C}:\Re(a\rho^n)<0\}$. Then
\begin{thm} \label{thm:P1:Intro:WellPosed}
The problem $\Pi(n,A,a,h,q_0)$ is well-posed if and only if $\eta_j(\rho)$ is entire and the ratio
\BE \label{eqn:P1:Intro:thm.WellPosed:Decay}
\frac{\eta_j(\rho)}{\DeltaP(\rho)}\to0 \mbox{ as }\rho\to\infty \mbox{ from within } D, \mbox{ away from the zeros of } \DeltaP.
\EE
for each $j$.
\end{thm}

We provide a small contribution to Fokas' method, making it fully algorithmic. We express the solution in terms of the PDE characteristic determinant, $\DeltaP$, the determinant of the matrix
\BE \label{eqn:P1:Intro:PDE.Characteristic.Matrix}
\M{\mathcal{A}}{k}{j}(\rho) = \begin{cases}\begin{array}{l}c_{(J_j-1)/2}(\rho)\left(\omega^{(n-1-[J_j-1]/2)(k-1)}\phantom{\displaystyle\sum_{r\in\widehat{J}^+}}\right. \\ \hspace{2.5mm} - \displaystyle\sum_{r\in\widehat{J}^+}\alpha_{\widehat{J}^+_r\hspace{0.5mm}(J_j-1)/2}\omega^{(n-1-r)(k-1)}(i\rho)^{(J_j-1)/2-r} \\ \hspace{2.5mm}\left.+ e^{-i\omega^{k-1}\rho}\displaystyle\sum_{r\in\widehat{J}^-}\alpha_{\widehat{J}^-_r\hspace{0.5mm}(J_j-1)/2}\omega^{(n-1-r)(k-1)}(i\rho)^{(J_j-1)/2-r}\right)\end{array} & J_j\mbox{ odd,} \\ \begin{array}{l}c_{J_j/2}(\rho)\left(-\omega^{(n-1-J_j/2)(k-1)}e^{-i\omega^{k-1}\rho}\phantom{\displaystyle\sum_{r\in\widehat{J}^+}}\right. \\ \hspace{2.5mm} - \displaystyle\sum_{r\in\widehat{J}^+}\beta_{\widehat{J}^+_r\hspace{0.5mm}J_j/2}\omega^{(n-1-r)(k-1)}(i\rho)^{J_j/2-r} \\ \hspace{2.5mm}\left.+ e^{-i\omega^{k-1}\rho}\displaystyle\sum_{r\in\widehat{J}^-}\beta_{\widehat{J}^-_r\hspace{0.5mm}J_j/2}\omega^{(n-1-r)(k-1)}(i\rho)^{J_j/2-r}\right)\end{array} & J_j\mbox{ even.}\end{cases}
\EE
The matrix $\mathcal{A}$ appears in the generalised spectral Dirichlet to Neumann map derived in Section~\ref{sec:P1:Implicit}. The application of the map to the formal result Theorem~\ref{thm:P1:Implicit:Formal} yields the following implicit equation for $q$, the solution of $\Pi$.

\begin{thm} \label{thm:P1:Intro:Implicit}
Let $\Pi(n,A,a,h,q_0)$ be well-posed with solution $q$. Then $q(x,t)$ may be expressed in terms of contour integrals of transforms of the boundary data, initial datum and solution at final time as follows:
\begin{multline} \label{thm:P1:Intro:thm.Implicit:q}
2\pi q(x,t)=\int_\mathbb{R}e^{i\rho x-a\rho^nt}\hat{q}_0(\rho)\d\rho - \int_{\partial D^+}e^{i\rho x-a\rho^nt}\sum_{j\in J^+}\frac{\zeta_j(\rho) - e^{a\rho^nT}\eta_j(\rho)}{\DeltaP(\rho)}\d\rho \\
- \int_{\partial D^-}e^{i\rho(x-1)-a\rho^nt}\sum_{j\in J^-}\frac{\zeta_j(\rho) - e^{a\rho^nT}\eta_j(\rho)}{\DeltaP(\rho)}\d\rho,
\end{multline}
where the sectors $D^\pm = D\cap\mathbb{C}^\pm$ and $D = \{\rho\in\mathbb{C}:\Re(a\rho^n)<0\}$.
\end{thm}

Equation~\eqref{thm:P1:Intro:thm.Implicit:q} gives only an implicit representation of the solution as the functions $\eta_j$ are defined in terms of the Fourier transform of the solution evaluated at final time, which is not a datum of the problem. Nevertheless the importance of the PDE characteristic determinant is clear. The integrands are meromorphic functions so $q$ depends upon their behaviour as $\rho\to\infty$ from within $D^\pm$ and upon their poles, which can only arise at zeros of $\DeltaP$. It is the behaviour at infinity that is used to characterise well-posedness in Theorem~\ref{thm:P1:Intro:WellPosed}, the proof of which is given in Section~\ref{sec:P1:WellPosed}.

In Section~\ref{sec:P1:Reps} we derive two representations of the solution of an initial-boundary value problem. Let $(\sigma_k)_{k\in\mathbb{N}}$ be a sequence containing each nonzero zero of $\DeltaP$ precisely once and define the index sets
\begin{gather*}
K^{\mathbb{R}}=\{k\in\mathbb{N}:\sigma_k\in\mathbb{R}\},\\
K^+=\{k\in\mathbb{N}:\Im\sigma_k\geq0\},\\
K^-=\{k\in\mathbb{N}:\Im\sigma_k<0\}.
\end{gather*}
Then the following theorems give representations of the solution to the problem $\Pi$.

\begin{thm} \label{thm:P1:Intro:Reps.Int}
Let the problem $\Pi(n,A,a,h,q_0)$ be well-posed. Then the solution $q$ may be expressed using contour integrals of transforms of the initial and boundary data by
\begin{multline} \label{eqn:P1:Intro:thm.Reps.Int:q}
q(x,t) = \frac{i}{2}\sum_{k\in K^+}\res_{\rho=\sigma_k}\frac{e^{i\rho x-a\rho^nt}}{\DeltaP(\rho)} \sum_{j\in J^+} \zeta_j(\rho) + \int_{\partial\widetilde{E}^+}e^{i\rho x-a\rho^nt}\sum_{j\in J^+}\frac{\zeta_j(\rho)}{\DeltaP(\rho)}\d\rho \\
+ \frac{i}{2}\sum_{k\in K^-}\res_{\rho=\sigma_k}\frac{e^{i\rho(x-1)-a\rho^nt}}{\DeltaP(\rho)} \sum_{j\in J^-} \zeta_j(\rho) + \int_{\partial\widetilde{E}^-}e^{i\rho(x-1)-a\rho^nt}\sum_{j\in J^-}\frac{\zeta_j(\rho)}{\DeltaP(\rho)}\d\rho \\
- \frac{1}{2\pi}\left\{\sum_{k\in K^\mathbb{R}} \int_{\Gamma_k} + \int_{\mathbb{R}} \right\} e^{i\rho x-a\rho^nt} \left( \frac{1}{\DeltaP(\rho)}-1 \right)H(\rho) \d\rho,
\end{multline}
\end{thm}

\begin{thm} \label{thm:P1:Intro:Reps.Ser}
Let $a=\pm i$ and let the problems $\Pi=\Pi(n,A,a,h,q_0)$ and $\Pi'=\Pi(n,A,-a,h,q_0)$ be well-posed. Then the solution $q$ of $\Pi$ may be expressed as a discrete series of transforms of the initial and boundary data by
\begin{multline} \label{eqn:P1:Intro:thm.Reps.Ser:q}
q(x,t) = \frac{i}{2}\sum_{k\in K^+}\res_{\rho=\sigma_k}\frac{e^{i\rho x-a\rho^nt}}{\DeltaP(\rho)} \sum_{j\in J^+} \zeta_j(\rho) \\ + \frac{i}{2}\sum_{k\in K^-}\res_{\rho=\sigma_k}\frac{e^{i\rho(x-1)-a\rho^nt}}{\DeltaP(\rho)} \sum_{j\in J^-} \zeta_j(\rho) \\
- \frac{1}{2\pi}\left\{\sum_{k\in K^\mathbb{R}} \int_{\Gamma_k} + \int_{\mathbb{R}} \right\} e^{i\rho x-a\rho^nt} \left( \frac{1}{\DeltaP(\rho)}-1 \right)H(\rho) \d\rho.
\end{multline}
\end{thm}

The final integral term in both equations~\eqref{eqn:P1:Intro:thm.Reps.Int:q} and~\eqref{eqn:P1:Intro:thm.Reps.Ser:q} depends upon $H$, a linear combination of $t$-transforms of the boundary data which evaluates to $0$ if $h=0$. Hence if $\Pi$ is a homogeneous initial-boundary value problem then the final term makes no contribution to equations~\eqref{eqn:P1:Intro:thm.Reps.Int:q} and~\eqref{eqn:P1:Intro:thm.Reps.Ser:q}.

Special cases of Theorem~\ref{thm:P1:Intro:Reps.Int} have appeared before but the representations differ from equation~\eqref{eqn:P1:Intro:thm.Reps.Int:q}. The result is shown for several specific examples in~\cite{FP2005a,Pel2005a}, including a second order problem with Robin boundary conditions. For simple boundary conditions, the result is mentioned in Remark~6 of~\cite{FP2001a} and Lemma~4.2 of~\cite{Pel2004a} contains the essence of the proof. Unlike earlier forms, equation~\eqref{eqn:P1:Intro:thm.Reps.Int:q} represents $q$ using discrete series as far as possible; only the parts of the integral terms that cannot be represented as series remain. This may not have any advantage for computation but is done to highlight the contrast with equation~\eqref{eqn:P1:Intro:thm.Reps.Ser:q}.

In Theorem~\ref{thm:P1:Intro:Reps.Ser} the well-posedness of $\Pi'$ is used to show that the first two integral terms of equation~\eqref{eqn:P1:Intro:thm.Reps.Int:q} evaluate to zero. Under the map $a\mapsto-a$, $D$ maps to $E$, the interior of its complement; we exploit this fact together with Theorem~\ref{thm:P1:Intro:WellPosed} to show the decay of
\BES
\frac{\zeta_j(\rho)}{\DeltaP(\rho)} \mbox{ as } \rho\to\infty \mbox{ from within } \widetilde{E}.
\EES
This maximally generalises of the arguments of Pelloni and Chilton in the sense that the deformation of contours cannot yield a series representation of the solution to $\Pi$ if $\Pi'$ is ill-posed.

Theorem~\ref{thm:P1:Intro:WellPosed} is useful because it reduces the complexity of the analysis necessary to prove that a particular initial-boundary value problem is well-posed but its use still requires some asymptotic analysis. It would be preferable to give a condition that may be validated by inspection of the boundary coefficient matrix and is sufficient for well-posedness. We discuss such criteria in Section~\ref{sec:P1:Alt}.

Section~\ref{sec:P1:Alt} also contains a proof of the following result, complementing Theorem~\ref{thm:P1:Intro:Reps.Ser}. This theorem highlights the essential difference between odd order problems, whose well-posedness depends upon the direction coefficient, and even order problems, whose well-posedness is determined by the boundary coefficient matrix only.

\begin{thm} \label{thm:P1:Intro:Even.Can.Deform.If.Well-posed}
Let $n$ be even and $a=\pm i$. Using the notation of Theorem~\ref{thm:P1:Intro:Reps.Ser}, the problem $\Pi'$ is well-posed if and only if $\Pi$ is well-posed.
\end{thm}

In Section~\ref{sec:P1:Spectrum} we investigate the \emph{PDE discrete spectrum}, the set of zeros of the PDE characteristic determinant. We prove a technical lemma, describing the distribution of the $\sigma_k$ which is used in the earlier sections. Under certain conditions we are able to exploit symmetry arguments to improve upon the general results Langer presents~\cite{Lan1931a} for the particular exponential polynomials of interest.

\section{Implicit solution of IBVP} \label{sec:P1:Implicit}

In Section~\ref{ssec:P1:Implicit:Fokas} we give the standard results of Fokas' unified transform method in the notation of this work. In Section~\ref{ssec:P1:Implicit:DtoN} we state and prove Lemma~\ref{lem:P1:Implicit:DtoNMap}, the generalised spectral Dirichlet to Neumann map. In Section~\ref{ssec:P1:Implicit:ApplyMap} we apply the map to the formal results of Section~\ref{ssec:P1:Implicit:Fokas}, concluding the proof of Theorem~\ref{thm:P1:Intro:Implicit}. The latter two sections contain formal definitions of many of the terms and much of the notation used throughout this work.

\subsection{Fokas' method} \label{ssec:P1:Implicit:Fokas}

The first steps of Fokas' transform method yield a formal representation for the solution of the initial-boundary value problem, given in the following

\begin{thm} \label{thm:P1:Implicit:Formal}
Let the initial-boundary value problem $\Pi(n,A,a,h,q_0)$ be well-posed. Then its solution $q$ may be expressed formally as the sum of three contour integrals,
\begin{multline} \label{eqn:P1:Implicit:thm.Formal:q}
q(x,t)=\frac{1}{2\pi}\left(\int_\mathbb{R}e^{i\rho x-a\rho^nt}\hat{q}_0(\rho)\d\rho - \int_{\partial D^+}e^{i\rho x-a\rho^nt}\sum_{j=0}^{n-1}c_j(\rho)\widetilde{f}_j(\rho)\d\rho\right. \\
\left.- \int_{\partial D^-}e^{i\rho(x-1)-a\rho^nt}\sum_{j=0}^{n-1}c_j(\rho)\widetilde{g}_j(\rho)\d\rho\right),
\end{multline}
where
\BE \label{eqn:P1:Implicit:thm.Formal:Definitions}
\begin{aligned}
\widetilde{f}_j(\rho)&=\int_0^Te^{a\rho^ns}f_j(s)\d s, & \widetilde{g}_j(\rho)&=\int_0^Te^{a\rho^ns}g_j(s)\d s, \\
f_j(t) &= \partial_x^jq(0,t), & g_j(t) &= \partial_x^jq(1,t), \\
\hat{q}_0(\rho)&=\int_0^1e^{-i\rho y}q_0(y)\d y, & c_j(\rho) &=-a\rho^n(i\rho)^{-(j+1)}.
\end{aligned}
\EE
\end{thm}

The above theorem is well established and its proof, via Lax pair and Riemann-Hilbert formalism, appears in~\cite{Fok2001a,Fok2008a,FP2001a}. We state it here without proof to highlight the difference in notation to previous publications. We use $\rho$ to denote the spectral parameter, in place of $k$ in the earlier work. We use $f_j$ and $g_j$ exclusively to denote the boundary functions; even for simple boundary conditions in which some of the boundary functions are equal to boundary data we denote the boundary data separately by $h_k$.

The transformed boundary functions are the $2n$ unknowns in equation~\eqref{eqn:P1:Implicit:thm.Formal:q}, of which at most $n$ may be explicitly specified by the boundary conditions~\eqref{eqn:P1:Intro:BC}. To determine the remaining $n$ or more we require a generalised Dirichlet to Neumann map in the form of Lemma~\ref{lem:P1:Implicit:DtoNMap}. This is derived from the boundary conditions and the global relation.

\begin{lem}[Global relation] \label{lem:P1:Implicit:GR}
Let $\Pi(n,A,a,h,q_0)$ be well-posed with solution $q$. Let
\BES
\hat{q}_T(\rho)=\int_0^1e^{-i\rho y}q(y,T)\d y
\EES
be the usual spatial Fourier transform of the solution evaluated at final time. Then the transformed functions $\hat{q}_0$, $\hat{q}_T$, $\widetilde{f}_j$ and $\widetilde{g}_j$ satisfy
\BE \label{lem:P1:Implicit:lem.GR:GR}
\sum_{j=0}^{n-1}c_j(\rho)\left( \tilde{f}_j(\rho)-e^{-i\rho}\tilde{g}_j(\rho) \right) = \hat{q}_0(\rho) - e^{a\rho^nT}\hat{q}_T(\rho),\qquad \rho\in\mathbb{C}.
\EE
\end{lem}

The global relation is derived using an application of Green's Theorem to the domain $[0,1]\times[0,T]$ in the aforementioned publications. As the $t$-transform,
\BE \label{eqn:P1:Implicit:tTransform}
\widetilde{X}(\rho)=\int_0^Te^{a\rho^nt}X(t)\d t,
\EE
is invariant under the map $\rho\mapsto \exp{(2j\pi i/n)}\rho$ for any integer $j$, the global relation provides a system of $n$ equations in the transformed functions to complement the boundary conditions.

\subsection{Generalised spectral Dirichlet to Neumann map} \label{ssec:P1:Implicit:DtoN}

We give a classification of boundary conditions and formally state the generalised spectral Dirichlet to Neumann map.

\begin{ntn} \label{ntn:P1:Implicit:Index.Sets}
Consider the problem $\Pi(n,A,a,h,q_0)$, which need not be well-posed. Define $\omega=\exp{(2\pi i/n)}$. Define the \emph{boundary coefficients} $\M{\alpha}{k}{j}$, $\M{\beta}{k}{j}$ to be the entries of $A$ such that
\BE \label{ntn:P1:Implicit:ntn.Index.Sets:Boundary.Coefficients}
\begin{pmatrix}\alpha_{1\hspace{0.5mm}n-1} & \beta_{1\hspace{0.5mm}n-1} & \alpha_{1\hspace{0.5mm}n-2} & \beta_{1\hspace{0.5mm}n-2}& \dots & \alpha_{1\hspace{0.5mm}0} & \beta_{1\hspace{0.5mm}0} \\ \alpha_{2\hspace{0.5mm}n-1} & \beta_{2\hspace{0.5mm}n-1} & \alpha_{2\hspace{0.5mm}n-2} & \beta_{2\hspace{0.5mm}n-2}& \dots & \alpha_{2\hspace{0.5mm}0} & \beta_{2\hspace{0.5mm}0} \\ \vdots & \vdots & \vdots & \vdots &  & \vdots & \vdots \\ \alpha_{n\hspace{0.5mm}n-1} & \beta_{n\hspace{0.5mm}n-1} & \alpha_{n\hspace{0.5mm}n-2} & \beta_{n\hspace{0.5mm}n-2}& \dots & \alpha_{n\hspace{0.5mm}0} & \beta_{n\hspace{0.5mm}0}\end{pmatrix}=A.
\EE
We define the following index sets and functions.

$\widehat{J}^+ = \{j\in\{0,1,\dots,n-1\}$ such that $\alpha_{k\hspace{0.5mm}j}$ is a pivot in $A$ for some $k\}$, the set of columns of $A$ relating to the left of the space interval which contain a pivot.

$\widehat{J}^- = \{j\in\{0,1,\dots,n-1\}$ such that $\beta_{k\hspace{0.5mm}j}$ is a pivot in $A$ for some $k\}$, the set of columns of $A$ relating to the right of the space interval which contain a pivot.

$\widetilde{J}^+ = \{0,1,\dots,n-1\}\setminus \widehat{J}^+$, the set of columns of $A$ relating to the left of the space interval which do not contain a pivot.

$\widetilde{J}^- = \{0,1,\dots,n-1\}\setminus \widehat{J}^-$, the set of columns of $A$ relating to the right of the space interval which do not contain a pivot.

$J = \{2j+1$ such that $j\in\widetilde{J}^+\}\cup\{2j$ such that $j\in\widetilde{J}^-\}$, an index set for the boundary functions whose corresponding columns in $A$ do not contain a pivot. Also, the decreasing sequence $(J_j)_{j=1}^n$ of elements of $J$.

$J' = \{2j+1$ such that $j\in\widehat{J}^+\}\cup\{2j$ such that $j\in\widehat{J}^-\} = \{0,1,\dots,2n-1\} \setminus J$, an index set for the boundary functions whose corresponding columns in $A$ contain a pivot. Also, the decreasing sequence $(J'_j)_{j=1}^n$ of elements of $J'$.

The functions $$V(\rho)=(V_1(\rho),V_2(\rho),\dots,V_n(\rho))^\T, \qquad V_j(\rho)=\begin{cases} \widetilde{f}_{(J_j-1)/2}(\rho)  & J_j \mbox{ odd,}  \\ \widetilde{g}_{J_j/2}(\rho)  & J_j \mbox{ even,}\end{cases}$$ the boundary functions whose corresponding columns in $A$ do not contain a pivot.

The functions $$W(\rho)=(W_1(\rho),W_2(\rho),\dots,W_n(\rho))^\T, \qquad W_j(\rho)=\begin{cases} \widetilde{f}_{(J'_j-1)/2}(\rho) & J'_j \mbox{ odd,} \\ \widetilde{g}_{J'_j/2}(\rho) & J'_j \mbox{ even,}\end{cases}$$ the boundary functions whose corresponding columns in $A$ contain a pivot.

$(\widehat{J}^+_j)_{j\in\widehat{J}^+}$, a sequence such that $\alpha_{\widehat{J}^+_j\hspace{0.5mm}j}$ is a pivot in $A$ when $j\in\widehat{J}^+$.

$(\widehat{J}^-_j)_{j\in\widehat{J}^-}$, a sequence such that $\beta_{\widehat{J}^-_j\hspace{0.5mm}j}$ is a pivot in $A$ when $j\in\widehat{J}^-$.
\end{ntn}

\begin{defn}[Classification of boundary conditions] \label{defn:P1:Implicit:BC.Classification}
The boundary conditions of the problem $\Pi(n,A,a,h,q_0)$ are said to be
\begin{enumerate}
\item{\emph{homogeneous} if $h=0$. Otherwise the boundary conditions are \emph{inhomogeneous}.}
\item{\emph{uncoupled} if
\begin{align*}
&\mbox{if } \M{\alpha}{k}{j} \mbox{ is a pivot in } A \mbox{ then } \M{\beta}{k}{r}=0 \hsforall r \mbox{ and} \\
&\mbox{if } \M{\beta}{k}{j} \mbox{ is a pivot in } A \mbox{ then } \M{\alpha}{k}{r}=0 \hsforall r.
\end{align*}
Otherwise we say that the boundary conditions are \emph{coupled}.}
\item{\emph{non-Robin} if
\BES
\hsforall k\in\{1,2,\dots,n\}, \mbox{ if } \M{\alpha}{k}{j}\neq0 \mbox{ or } \M{\beta}{k}{j}\neq0 \mbox{ then } \M{\alpha}{k}{r}=\M{\beta}{k}{r}=0\hsforall r\neq j,
\EES
that is each contains only one order of partial derivative. Otherwise we say that boundary condition is \emph{of Robin type}. Note that whether boundary conditions are of Robin type or not is independent of whether they are coupled, unlike Duff's definition~\cite{Duf1956a}.}
\item{\emph{simple} if they are uncoupled and non-Robin.}
\end{enumerate}
\end{defn}

The terms `generalised' and `spectral' are prefixed to the name `Dirichlet to Neumann map' of the Lemma below to avoid confusion regarding its function.

{\bfseries Generalised:} The boundary conditions we study are considerably more complex than those considered in~\cite{Chi2006a,Fok2001a,FP2001a,FP2005a,Pel2004a,Pel2005a}. Indeed, as $A$ may specify any linear boundary conditions, the known boundary functions may not be `Dirichlet' (zero order) and the unknown boundary functions need not be `Neumann' (first order). Further, if $A$ has more than $n$ non-zero entries then the lemma must be capable of expressing more than $n$ unknown boundary functions in terms of fewer than $n$ known boundary data.

{\bfseries Spectral:} Owing to the form of equation~\eqref{eqn:P1:Implicit:thm.Formal:q} we are interested not in the boundary functions themselves but in their $t$-transforms, as defined in equations~\eqref{eqn:P1:Implicit:thm.Formal:Definitions}. It is possible, though unnecessarily complicated, to perform a generalized Dirichlet to Neumann map in real time and subsequently transform to spectral time but, as the global relation is in spectral time, to do so requires the use of an inverse spectral transform. Instead, we exploit the linearity of the $t$-transform~\eqref{eqn:P1:Implicit:tTransform}, applying it to the boundary conditions, and derive the map in spectral time.

The crucial component of the lemma is given in the following

\begin{defn} \label{defn:P1:Implicit:DeltaP}
Let $\Pi(n,A,a,h,q_0)$ be an initial-boundary value problem having the properties $(\Pi1)$--$(\Pi4)$ but not necessarily well-posed. We define the \emph{PDE characteristic matrix} $\mathcal{A}(\rho)$ by equation~\eqref{eqn:P1:Intro:PDE.Characteristic.Matrix} and the \emph{PDE characteristic determinant} to be the entire function
\BE \label{eqn:P1:Implicit:defn.DeltaP:DeltaP}
\DeltaP(\rho)=\det\mathcal{A}(\rho).
\EE
\end{defn}

\begin{lem}[Generalised spectral Dirichlet to Neumann map] \label{lem:P1:Implicit:DtoNMap}
Let $\Pi(n,A,a,h,q_0)$ be well-posed with solution $q$. Then
\begin{enumerate}
\item{The vector $V$ of transformed boundary functions satisfies the \emph{reduced global relation}
\BE \label{eqn:P1:Implicit:lem.DtoNMap:RGR}
\mathcal{A}(\rho)V(\rho) = U(\rho) - e^{a\rho^nT}\begin{pmatrix}\hat{q}_T(\rho)\\\vdots\\\hat{q}_T(\omega^{n-1}\rho)\end{pmatrix},
\EE
where
\begin{align} \label{eqn:P1:Implicit:lem.DtoNMap:U}
U(\rho) &= (u(\rho,1),u(\rho,2),\dots,u(\rho,n))^\T, \\ \label{eqn:P1:Implicit:lem.DtoNMap:u}
u(\rho,k) &= \hat{q}_0(\omega^{k-1}\rho) - \sum_{l\in\widehat{J}^+} c_l(\omega^{k-1}\rho)\widetilde{h}_{\widehat{J}^+_l}(\rho) + e^{-i\omega^{k-1}\rho}\sum_{l\in\widehat{J}^-} c_l(\omega^{k-1}\rho)\widetilde{h}_{\widehat{J}^-_l}(\rho)
\end{align}
and $\widetilde{h}_j$ is the function obtained by applying the $t$-transform~\eqref{eqn:P1:Implicit:tTransform} to the boundary datum $h_j$.}
\item{The PDE characteristic matrix is of full rank, is independent of $h$ and $q_0$ and differing values of $a$ only scale $\mathcal{A}$ by a nonzero constant factor.}
\item{The vectors $V$ and $W$ of transformed boundary functions satisfy the \emph{reduced boundary conditions}
\BE \label{eqn:P1:Implicit:lem.DtoNMap:RBC}
W(\rho) = \left(\widetilde{h}_1(\rho),\widetilde{h}_2(\rho),\dots,\widetilde{h}_n(\rho)\right)^\T - \widehat{A} V(\rho),
\EE
where the \emph{reduced boundary coefficient matrix} is given by
\BE \label{eqn:P1:Implicit:lem.DtoNMap:RBC.Matrix}
\widehat{A}_{k\hspace{0.5mm}j} = \begin{cases} \M{\alpha}{k}{(J_j-1)/2} & J_j \mbox{ odd,} \\ \M{\beta}{k}{J_j/2} & J_j \mbox{ even.}\end{cases}
\EE}
\end{enumerate}
\end{lem}

\begin{proof}
Applying the $t$-transform~\eqref{eqn:P1:Implicit:tTransform} to each line of the boundary conditions~\eqref{eqn:P1:Intro:BC} yields a system of $n$ equations in the transformed boundary functions. As $A$ is in reduced row-echelon form it is possible to split the vector containing all of the transformed boundary functions into the two vectors $V$ and $W$, justifying the reduced boundary conditions.

The reduced boundary conditions may also be written
\begin{align} \label{eqn:GettingA:General:Main.Lemma.Proof.Reduced.GR1}
\widetilde{f}_j(\rho) &= \widetilde{h}_{\widehat{J}^+_j}(\rho) - \sum_{r\in\widetilde{J}^+}\alpha_{\widehat{J}^+_j\hspace{0.5mm}r}\widetilde{f}_r(\rho) - \sum_{r\in\widetilde{J}^-}\beta_{\widehat{J}^+_j\hspace{0.5mm}r}\widetilde{g}_r(\rho), & \mbox{for } j&\in\widehat{J}^+ \mbox{ and} \\ \label{eqn:GettingA:General:Main.Lemma.Proof.Reduced.GR2}
\widetilde{g}_j(\rho) &= \widetilde{h}_{\widehat{J}^-_j}(\rho) - \sum_{r\in\widetilde{J}^+}\alpha_{\widehat{J}^-_j\hspace{0.5mm}r}\widetilde{f}_r(\rho) - \sum_{r\in\widetilde{J}^-}\beta_{\widehat{J}^-_j\hspace{0.5mm}r}\widetilde{g}_r(\rho), & \mbox{for } j&\in\widehat{J}^-.
\end{align}
As the $t$-transform is invariant under the map $\rho\mapsto\omega^j\rho$ for any integer $j$, the global relation Lemma~\ref{lem:P1:Implicit:GR} yields the system
\BES
\sum_{j=0}^{n-1}c_{j}(\rho)\omega^{(n-1-j)r}\widetilde{f}_{j}(\rho) - \sum_{j=0}^{n-1}e^{-i\omega^r\rho}c_{j}(\rho)\omega^{(n-1-j)r}\widetilde{g}_{j}(\rho) = \hat{q}_0(\omega^r\rho) - e^{a\rho^nT}\hat{q}_T(\omega^r\rho),
\EES
for $r\in\{0,1,\dots,n-1\}$. Using the fact $\widehat{J}^+\cup\widetilde{J}^+=\widehat{J}^-\cup\widetilde{J}^-=\{0,1,\dots,n-1\}$ we split the sums on the left hand side to give
\begin{multline*}
\sum_{j\in\widehat{J}^+}c_{j}(\rho)\omega^{(n-1-j)r}\widetilde{f}_{j}(\rho) + \sum_{j\in\widetilde{J}^+}c_{j}(\rho)\omega^{(n-1-j)r}\widetilde{f}_{j}(\rho) \\
- \sum_{j\in\widehat{J}^-}e^{-i\omega^r\rho}c_{j}(\rho)\omega^{(n-1-j)r}\widetilde{g}_{j}(\rho) - \sum_{j\in\widetilde{J}^-}e^{-i\omega^r\rho}c_{j}(\rho)\omega^{(n-1-j)r}\widetilde{g}_{j}(\rho) \\
= \hat{q}_0(\omega^r\rho) - e^{a\rho^nT}\hat{q}_T(\omega^r\rho),
\end{multline*}
for $r\in\{0,1,\dots,n-1\}$. Substituting equations~\eqref{eqn:GettingA:General:Main.Lemma.Proof.Reduced.GR1} and~\eqref{eqn:GettingA:General:Main.Lemma.Proof.Reduced.GR2} and interchanging the summations we obtain the reduced global relation.

The latter statement of (ii) is a trivial observation from the form of the PDE characteristic matrix. A full proof that $\mathcal{A}$ is full rank is given in the proof of Lemma~2.17 of~\cite{Smi2011a}.
\end{proof}

\subsection{Applying the map} \label{ssec:P1:Implicit:ApplyMap}

We solve the system of linear equations~\eqref{eqn:P1:Implicit:lem.DtoNMap:RGR} for $V$ using Cramer's rule hence, by equation~\eqref{eqn:P1:Implicit:lem.DtoNMap:RBC}, determining $W$ also.

\begin{ntn} \label{ntn:P1:Implicit:DeltaP}
Denote by $\widehat{\zeta}_j(\rho)$ the determinant of the matrix obtained by replacing the $j$\Th column of the PDE characteristic matrix with the vector $U(\rho)$ and denote by $\widehat{\eta}_j(\rho)$ the determinant of the matrix obtained by replacing the $j$\Th column of the PDE characteristic matrix with the vector $(\hat{q}_T(\rho),\hat{q}_T(\omega\rho),\dots,\hat{q}_T(\omega^{n-1}\rho))^\T$ for $j\in\{1,2,\dots,n\}$ and $\rho\in\mathbb{C}$. Define
\BE \label{eqn:P1:Implicit:defn.DeltaP:zetajhat.etajhat} \begin{split}
\widehat{\zeta}_j(\rho) &= \widetilde{h}_{j-n}(\rho) - \sum_{k=1}^{n}\M{\widehat{A}}{j-n}{k}\widehat{\zeta}_{k}(\rho), \\
\widehat{\eta}_j(\rho)  &= \widetilde{h}_{j-n}(\rho) - \sum_{k=1}^{n}\M{\widehat{A}}{j-n}{k}\widehat{\eta}_{k}(\rho),
\end{split}\EE
for $j\in\{n+1,n+2,\dots,2n\}$ and $\rho\in\mathbb{C}$. Define
\BE \label{eqn:P1:Implicit:defn.DeltaP:zetaj.etaj}
\zeta_j(\rho) = \begin{cases}c_{(J_j-1)/2}(\rho)\widehat{\zeta}_j(\rho) \\ c_{J_j/2}(\rho)\widehat{\zeta}_j(\rho) \\ c_{(J'_{j-n}-1)/2}(\rho)\widehat{\zeta}_j(\rho) \\ c_{J'_{j-n}/2}(\rho)\widehat{\zeta}_j(\rho)\end{cases} \qquad
\eta_j(\rho) = \begin{cases}c_{(J_j-1)/2}(\rho)\widehat{\eta}_j(\rho) & J_j \mbox{ odd,} \\ c_{J_j/2}(\rho)\widehat{\eta}_j(\rho) & J_j \mbox{ even,} \\ c_{(J'_{j-n}-1)/2}(\rho)\widehat{\eta}_j(\rho) & J'_{j-n} \mbox{ odd,} \\ c_{J'_{j-n}/2}(\rho)\widehat{\eta}_j(\rho) & J'_{j-n} \mbox{ even,}\end{cases}
\EE
for $\rho\in\mathbb{C}$ and define the index sets
\begin{align*}
J^+ &= \{j:J_j \mbox{ odd}\} \cup\{n+j:J'_j \mbox{ odd} \}, \\
J^- &= \{j:J_j \mbox{ even}\}\cup\{n+j:J'_j \mbox{ even}\}.
\end{align*}
\end{ntn}

The generalised spectral Dirichlet to Neumann map Lemma~\ref{lem:P1:Implicit:DtoNMap} and Cramer's rule yield expressions for the transformed boundary functions:
\BE \label{eqn:P1:Implicit:Zeta.Eta.f.g}
\frac{\zeta_j(\rho) - e^{a\rho^nT}\eta_j(\rho)}{\DeltaP(\rho)} = \begin{cases}c_{(J_j-1)/2}(\rho)\widetilde{f}_{(J_j-1)/2}(\rho) & J_j \mbox{ odd,} \\ c_{J_j/2}(\rho)\widetilde{g}_{J_j/2}(\rho) & J_j \mbox{ even,} \\ c_{(J'_{j-n}-1)/2}(\rho)\widetilde{f}_{(J'_{j-n}-1)/2}(\rho) & J'_{j-n} \mbox{ odd,} \\ c_{J'_{j-n}/2}(\rho)\widetilde{g}_{J'_{j-n}/2}(\rho) & J'_{j-n} \mbox{ even,}\end{cases}
\EE
hence
\begin{align*}
\sum_{j=0}^{n-1}c_j(\rho)\widetilde{f}_j(\rho) &= \sum_{j\in J^+}\frac{\zeta_j(\rho) - e^{a\rho^nT}\eta_j(\rho)}{\DeltaP(\rho)}, \\
\sum_{j=0}^{n-1}c_j(\rho)\widetilde{g}_j(\rho) &= \sum_{j\in J^-}\frac{\zeta_j(\rho) - e^{a\rho^nT}\eta_j(\rho)}{\DeltaP(\rho)}.
\end{align*}
Substituting these equations into Theorem~\ref{thm:P1:Implicit:Formal} completes the proof of Theorem~\ref{thm:P1:Intro:Implicit}.

\begin{rmk}
There are several simplifications of the above definitions for specific types of boundary conditions.

If the boundary conditions are simple, as studied in~\cite{Pel2004a}, then $\widehat{A}=0$. Hence, if the boundary conditions are simple and homogeneous then $\zeta_j=\eta_j=0$ for each $j>n$.

Non-Robin boundary conditions admit a significantly simplified form of the PDE characteristic matrix; see equation~(2.2.5) of~\cite{Smi2011a}.

For homogeneous boundary conditions, $\eta_j$ is $\zeta_j$ with $\hat{q}_T$ replacing $\hat{q}_0$.
\end{rmk}

\begin{rmk} \label{rmk:P1:Implicit:General.PDE}
It is possible to extend the results above to initial-boundary value problems for a more general linear, constant-coefficient evolution equation,
\BE \label{eqn:P1:Implicit:rmk.General.PDE:PDE}
\partial_tq(x,t) + \sum_{j=0}^na_j(-i\partial_x)^jq(x,t) = 0,
\EE
with leading coefficient $a_n$ having the properties of $a$. In this case the spectral transforms must be redefined with $\sum_{j=0}^na_j\rho^j$ replacing $a\rho^n$ and the form of the boundary coefficient matrix also changes. The $\omega^X$ appearing in equation~\eqref{eqn:P1:Intro:PDE.Characteristic.Matrix} represent a rotation by $2X\pi/n$, corresponding to a map between simply connected components of $D$. The partial differential equation~\eqref{eqn:P1:Implicit:rmk.General.PDE:PDE} has dispersion relation $\sum_{j=0}^na_j\rho^n$ so $D$ is not simply a union of sectors but a union of sets that are asymptotically sectors; see Lemma~1.1 of~\cite{FS1999a}. Hence we replace $\omega^X$ with a biholomorphic map between the components of $D$.
\end{rmk}

\section{New characterisation of well-posedness} \label{sec:P1:WellPosed}

This section provides a proof of Theorem~\ref{thm:P1:Intro:WellPosed}. The first subsection justifies that the decay condition is satisfied by all well-posed problems. The second subsection proves that the decay condition is sufficient for well-posedness.

We clarify the definitions of $\widetilde{D}$ and $\widetilde{E}$ from Section~\ref{sec:P1:Intro}. By Lemma~\ref{lem:P1:Spectrum:Properties}, there exists some $\varepsilon>0$ such that the pairwise intersection of closed discs of radius $\varepsilon$ centred at zeros of $\DeltaP$ is empty. We define
\BES
\widetilde{D}=D\setminus\bigcup_{k\in\mathbb{N}}\overline{B}(\sigma_k,\varepsilon), \quad \widetilde{E}=E\setminus\bigcup_{k\in\mathbb{N}}\overline{B}(\sigma_k,\varepsilon).
\EES

\subsection{Well-posedness $\Rightarrow$ decay} \label{ssec:P1:WellPosed:WP.Implies.Decay}

As the problem is well-posed, the solution evaluated at final time $q_T\in C^\infty[0,1]$ hence $\hat{q}_T$ and $\eta_j$ are entire. Similarly, $f_k,g_k\in C^\infty[0,T]$ hence $\widetilde{f}_k,\widetilde{g}_k$ are entire and decay as $\rho\to\infty$ from within $D$. Hence, by equation~\eqref{eqn:P1:Implicit:Zeta.Eta.f.g},
\BE \label{eqn:P1:WellPosed:WP.Implies.Decay.a}
\frac{\zeta_j(\rho) - e^{a\rho^nT}\eta_j(\rho)}{\DeltaP(\rho)c_k(\rho)}
\EE
is entire and decays as $\rho\to0$ from within $D$ for each $j\in\{1,2,\dots,2n\}$, where $k$ depends upon $j$.

We define the new complex set
\BES
\mathcal{D} = \{\rho\in D \mbox{ such that }-\Re(a\rho^nT)>2n|\rho|\}.
\EES
As $\mathcal{D}\subset D$, the ratio~\eqref{eqn:P1:WellPosed:WP.Implies.Decay.a} is analytic on $\mathcal{D}$ and decays as $\rho\to\infty$ from within $\mathcal{D}$. For $p\in\{1,2,\dots,n\}$, let $D_p$ be the $p$\Th simply connected component of $D$ encountered when moving anticlockwise from the positive real axis and let $\widetilde{D}_p=\widetilde{D}\cap D_p$. Then for each $p\in\{1,2,\dots,n\}$ there exists $R>0$ such that the set
\BES
\mathcal{D}_p = \left(\widetilde{D}_p\cap \mathcal{D}\right)\setminus\overline{B}(0,R)
\EES
is simply connected, open and unbounded.

By definition, $\DeltaP(\rho)$ is an exponential polynomial whose terms are each
\BES
W(\rho)e^{-i\sum_{y\in Y}\omega^r\rho}
\EES
where $W$ is a monomial of degree at least $1$ and $Y\subset\{0,1,2,\dots,n-1\}$ is an index set. Hence
\BES
\frac{1}{\DeltaP(\rho)} = o(e^{n|\rho|}\rho^{-1}) \mbox{ as } \rho\to\infty \mbox{ or as } \rho\to0.
\EES
As $\zeta_j$ and $\eta_j$ also grow no faster than $o(e^{n|\rho|})$, the ratios
\BES
\frac{\zeta_j(\rho)}{\DeltaP(\rho)c_k(\rho)}, \quad \frac{\eta_j(\rho)}{\DeltaP(\rho)c_k(\rho)} = o(e^{2n|\rho|}\rho^{-1}), \mbox{ as } \rho\to\infty.
\EES
Hence the ratio
\BE \label{eqn:P1:WellPosed:WP.Implies.Decay.d}
\frac{e^{a\rho^nT}\eta_j(\rho)}{\DeltaP(\rho)c_k(\rho)}
\EE
decays as $\rho\to\infty$ from within $\mathcal{D}$ and away from the zeros of $\DeltaP$. However the ratio
\BE \label{eqn:P1:WellPosed:WP.Implies.Decay.b}
\frac{\zeta_j(\rho)}{\DeltaP(\rho)c_k(\rho)}
\EE
is the sum of ratios~\eqref{eqn:P1:WellPosed:WP.Implies.Decay.a} and~\eqref{eqn:P1:WellPosed:WP.Implies.Decay.d} hence it also decays as $\rho\to\infty$ from within $D'$ and away from the zeros of $\DeltaP$.

The terms in each of $\zeta_j(\rho)$ and $\DeltaP(\rho)$ are exponentials, each of which either decays or grows as $\rho\to\infty$ from within one of the simply connected components $\widetilde{D}_p$ of $\widetilde{D}$. Hence as $\rho\to\infty$ from within a particular component $\widetilde{D}_p$ the ratio~\eqref{eqn:P1:WellPosed:WP.Implies.Decay.b} either decays or grows. But, as observed above, these ratios all decay as $\rho\to\infty$ from within each $\mathcal{D}_p$. Hence the ratio~\eqref{eqn:P1:WellPosed:WP.Implies.Decay.b} decays as $\rho\to\infty$ from within $\widetilde{D}_p$.

Now it is a simple observation that the ratio
\BE \label{eqn:P1:WellPosed:WP.Implies.Decay.c}
\frac{\eta_j(\rho)}{\DeltaP(\rho)c_k(\rho)}
\EE
must also decay as $\rho\to\infty$. Indeed ratio~\eqref{eqn:P1:WellPosed:WP.Implies.Decay.c} is the same as ratio~\eqref{eqn:P1:WellPosed:WP.Implies.Decay.b} but with $\hat{q}_T(\omega^{k-1}\rho)$ replacing $u(\rho,k)$ and, as observed above, $q_T\in C^\infty[0,1]$ also. Finally, the exponentials in $\eta_j$ and $\DeltaP$ ensure that the ratio
\BE \label{eqn:P1:WellPosed:WP.Implies.Decay.e}
\frac{\eta_j(\rho)}{\DeltaP(\rho)}
\EE
also decays as $\rho\to\infty$ from within $\widetilde{D}_p$. Indeed the transforms that multiply each term in $\eta_j$ ensure that the decay of ratio~\eqref{eqn:P1:WellPosed:WP.Implies.Decay.c} must come from the decay of ratio~\eqref{eqn:P1:WellPosed:WP.Implies.Decay.e}, not from $1/c_k(\rho)$.

\subsection{Decay $\Rightarrow$ well-posedness} \label{ssec:P1:WellPosed:Decay.Implies.WP}

Many of the definitions of Section~\ref{sec:P1:Implicit} require the problem $\Pi(n,A,a,h,q_0)$ to be well-posed. The statement of the following Lemma clarifies what is meant by $\eta_j$ when $\Pi$ is not known to be well-posed a priori and the result is the principal tool in the proof of Theorem~\ref{thm:P1:Intro:WellPosed}.

\begin{lem} \label{lem:P1:WellPosed:Ass.D.Implies.Admissible}
Consider the problem $\Pi(n,A,a,h,q_0)$ with associated PDE characteristic matrix $\mathcal{A}$ whose determinant is $\DeltaP$. Let the polynomials $c_j$ be defined by $c_j(\rho) = -a\rho^n(i\rho)^{-(j+1)}$. Let $U:\mathbb{C}\to\mathbb{C}$ be defined by equation~\eqref{eqn:P1:Implicit:lem.DtoNMap:U} and let $\widehat{A}\in\mathbb{R}^{n\times n}$ be defined by equation~\eqref{eqn:P1:Implicit:lem.DtoNMap:RBC}. Let $\zeta_j,\eta_j:\mathbb{C}\to\mathbb{C}$ be defined by Notation~\ref{ntn:P1:Implicit:DeltaP}, where $q_T:[0,1]\to\mathbb{C}$ is some function such that $\eta_j$ is entire and the decay condition~\eqref{eqn:P1:Intro:thm.WellPosed:Decay} is satisfied. Let the functions $\widetilde{f}_j,\widetilde{g}_j:\mathbb{C}\to\mathbb{C}$ be defined by equation~\eqref{eqn:P1:Implicit:Zeta.Eta.f.g}. Let $f_j,g_j:[0,T]\to\mathbb{C}$ be the functions for which
\BE \label{eqn:P1:WellPosed:lem.Ass.D.Implies.Admissible:Defn.fj}
\widetilde{f}_j(\rho) = \int_0^Te^{a\rho^nt}f_j(t)\d t,\quad \widetilde{g}_j(\rho) = \int_0^Te^{a\rho^nt}g_j(t)\d t, \quad \rho\in\mathbb{C}.
\EE
Then $\{f_j,g_j:j\in\{0,1,\dots,n-1\}\}$ is an admissible set in the sense of Definition~1.3 of~\cite{FP2001a}.
\end{lem}

\begin{proof}
By equation~\eqref{eqn:P1:Implicit:Zeta.Eta.f.g} and the definition of the index sets $J^\pm$ in Notation~\ref{ntn:P1:Implicit:DeltaP} we may write equations~(1.13) and~(1.14) of~\cite{FP2001a} as
\begin{align} \label{eqn:App:Ass.D.Implies.Admissible.Lem:Alt.F}
\widetilde{F}(\rho) &= \sum_{j\in J^+}{\frac{\zeta_j(\rho)-e^{a\rho^nT}\eta_j(\rho)}{\DeltaP(\rho)}}, \\ \label{eqn:App:Ass.D.Implies.Admissible.Lem:Alt.G}
\widetilde{G}(\rho) &= \sum_{j\in J^-}{\frac{\zeta_j(\rho)-e^{a\rho^nT}\eta_j(\rho)}{\DeltaP(\rho)}}.
\end{align}
By Cramer's rule and the calculations in the proof of Lemma~\ref{lem:P1:Implicit:DtoNMap}, equation~(1.17) of~\cite{FP2001a} is satisfied.

As $\eta_j$ is entire, $\hat{q}_T$ is entire so, by the standard results on the inverse Fourier transform, $q_T:[0,1]\to\mathbb{C}$, defined by
\BES
q_T(x) = \frac{1}{2\pi}\int_{\mathbb{R}}e^{i\rho x}\hat{q}_T(\rho)\d \rho,
\EES
is a $C^\infty$ smooth function.

We know $\zeta_j$ is entire by construction and $\eta_j$ is entire by assumption hence $\widetilde{F}$ and $\widetilde{G}$ are meromorphic on $\mathbb{C}$ and analytic on $\widetilde{D}$. By the definition of $D$ and the decay assumption
\BES
\frac{e^{a\rho^nT}\eta_j(\rho)}{\DeltaP(\rho)}\to0 \mbox{ as } \rho\to\infty \mbox{ from within } \widetilde{D}.
\EES
As $\hat{q}_0$ and $\widetilde{h}_j$ are entire so is $U$. As $\hat{q}_T$ is also entire and the definitions of $\zeta_j$ and $\eta_j$ differ only by which of these functions appears, the ratio $\zeta_k(\rho)/\DeltaP(\rho)\to0$ as $\rho\to\infty$ from within $\widetilde{D}$ also. This establishes that
\BES
\frac{\zeta_j(\rho)-e^{a\rho^nT}\eta_j(\rho)}{\DeltaP(\rho)}\to0 \mbox{ as } \rho\to\infty \mbox{ from within } \widetilde{D}.
\EES
Hence, by equations~\eqref{eqn:App:Ass.D.Implies.Admissible.Lem:Alt.F} and~\eqref{eqn:App:Ass.D.Implies.Admissible.Lem:Alt.G}, $\widetilde{F}(\rho),\widetilde{G}(\rho)\to0$ as $\rho\to\infty$ within $\widetilde{D}$.

An argument similar to that in Example~7.4.6 of~\cite{AF1997a} yields
\begin{align*}
f_j(t) &= -\frac{i^j}{2\pi}\int_{\partial D}\rho^je^{-a\rho^nt}\widetilde{F}(\rho)\d \rho, \\
g_j(t) &= -\frac{i^j}{2\pi}\int_{\partial D}\rho^je^{-a\rho^nt}\widetilde{G}(\rho)\d \rho.
\end{align*}
Because $\widetilde{F}(\rho),\widetilde{G}(\rho)\to0$ as $\rho\to\infty$ within $\widetilde{D}$, these definitions guarantee that $f_j$ and $g_j$ are $C^\infty$ smooth.

The compatibility of the $f_j$ and $g_j$ with $q_0$ is ensured by the compatibility condition $(\Pi4)$.
\end{proof}

The desired result is now a restatement of Theorems~1.1 and~1.2 of~\cite{FP2001a}. For this reason we refer the reader to the proof presented in Section~4 of that publication. The only difference is that we make use of Lemma~\ref{lem:P1:WellPosed:Ass.D.Implies.Admissible} in place of Proposition~4.1.

\section{Representations of the solution} \label{sec:P1:Reps}

The proofs of Theorems~\ref{thm:P1:Intro:Reps.Int} and~\ref{thm:P1:Intro:Reps.Ser} are similar calculations. In Section~\ref{ssec:P1:Reps:Ser} we present the derivation of the series representation and, in Section~\ref{ssec:P1:Reps:Int}, note the way this argument may be adapted to yield the integral representation. We derive the result in the case $n$ odd, $a=i$; the other cases are almost identical.

\subsection{Series Representation} \label{ssec:P1:Reps:Ser}

As $\Pi$ is well-posed, Theorem~\ref{thm:P1:Intro:Implicit} holds. We split the latter two integrals of equation~\eqref{thm:P1:Intro:thm.Implicit:q} into parts whose integrands contain the data, that is $\zeta_j$, and parts whose integrands contain the solution evaluated at final time, that is $\eta_j$.
\begin{multline} \label{eqn:P1:Reps:Ser:q.Implicit}
2\pi q(x,t) = \int_\mathbb{R}e^{i\rho x-i\rho^nt}\hat{q}_0(\rho)\d\rho + \left\{\int_{\partial E^+}-\int_{\mathbb{R}}\right\}e^{i\rho x-i\rho^nt}\sum_{j\in J^+}\frac{\zeta_j(\rho)}{\DeltaP(\rho)}\d\rho \\
+ \int_{\partial D^+}e^{i\rho x+i\rho^n(T-t)}\sum_{j\in J^+}\frac{\eta_j(\rho)}{\DeltaP(\rho)}\d\rho + \left\{\int_{\partial E^-}+\int_{\mathbb{R}}\right\}e^{i\rho(x-1)-i\rho^nt}\sum_{j\in J^-}\frac{\zeta_j(\rho)}{\DeltaP(\rho)}\d\rho \\
+ \int_{\partial D^-}e^{i\rho(x-1)+i\rho^n(T-t)}\sum_{j\in J^-}\frac{\eta_j(\rho)}{\DeltaP(\rho)}\d\rho.
\end{multline}

As $\Pi'$ is well-posed, Theorem~\ref{thm:P1:Intro:WellPosed} ensures the ratios
\BES
\frac{\eta'_j(\rho)}{\DeltaPprime(\rho)}\to0 \mbox{ as }\rho\to\infty \mbox{ from within } \widetilde{D}',
\EES
for each $j$. By definition $E=D'$ and, by statement (ii) of Lemma~\ref{lem:P1:Implicit:DtoNMap}, the zeros of $\mathcal{A}'$ are precisely the zeros of $\mathcal{A}$ hence $\widetilde{E}=\widetilde{D}'$. Define $\xi_j(\rho)$ to be the function obtained by replacing $\hat{q}'_T(\omega^{k-1}\rho)$ with $u(\rho,k)$ in the definition of $\eta'_j(\rho)$. As $q'_T$, $q_0$ and $h_j$ are all smooth functions, $\xi_j$ has precisely the same decay properties of $\eta'_j$. But $\xi_j=\zeta_j$ by definition. Hence the well-posedness of $\Pi'$ is equivalent to
\BE \label{eqn:P1:Reps:Ser:zetaj.Decay}
\frac{\zeta_j(\rho)}{\DeltaP(\rho)}\to0 \mbox{ as }\rho\to\infty \mbox{ from within } \widetilde{E},
\EE
for each $j$. The decay property obtained by applying Theorem~\ref{thm:P1:Intro:WellPosed} directly to $\Pi$ together with the decay property~\eqref{eqn:P1:Reps:Ser:zetaj.Decay} permits the use of Jordan's Lemma to deform the contours of integration over $\widetilde{D}^\pm$ and $\widetilde{E}^\pm$ in equation~\eqref{eqn:P1:Reps:Ser:q.Implicit} to obtain
\begin{multline} \label{eqn:P1:Reps:Ser:q.Integrals.at.zeros}
2\pi q(x,t) = \int_\mathbb{R}e^{i\rho x-i\rho^nt}\hat{q}_0(\rho)\d\rho + \left\{\int_{\partial (E^+\setminus\widetilde{E}^+)}-\int_{\mathbb{R}}\right\}e^{i\rho x-i\rho^nt}\sum_{j\in J^+}\frac{\zeta_j(\rho)}{\DeltaP(\rho)}\d\rho \\
+ \int_{\partial (D^+\setminus\widetilde{D}^+)}e^{i\rho x+i\rho^n(T-t)}\sum_{j\in J^+}\frac{\eta_j(\rho)}{\DeltaP(\rho)}\d\rho \\
+ \left\{\int_{\partial (E^-\setminus\widetilde{E}^-)}+\int_{\mathbb{R}}\right\}e^{i\rho(x-1)-i\rho^nt}\sum_{j\in J^-}\frac{\zeta_j(\rho)}{\DeltaP(\rho)}\d\rho \\
+ \int_{\partial (D^-\setminus\widetilde{D}^-)}e^{i\rho(x-1)+i\rho^n(T-t)}\sum_{j\in J^-}\frac{\eta_j(\rho)}{\DeltaP(\rho)}\d\rho.
\end{multline}
Indeed, $\zeta_j$, $\eta_j$ and $\DeltaP$ are entire functions hence the ratios can have poles only at the zeros of $\DeltaP$, neighbourhoods of which are excluded from $\widetilde{D}^\pm$ and $\widetilde{E}^\pm$ by definition. Finally, the exponential functions in the integrands each decay as $\rho\to\infty$ from within the sectors enclosed by their respective contour of integration.

The right hand side of equation~\eqref{eqn:P1:Reps:Ser:q.Integrals.at.zeros} is the sum of three integrals over $\mathbb{R}$ and four others. The former may be combined into a single integral using the following lemma, whose proof appears at the end of this section.
\begin{lem} \label{lem:P1:Reps:q0}
Let $\Pi(n,A,a,q_0,h)$ be well-posed. Then
\BE \label{eqn:P1:Repd:lem.q0:Statement}
\sum_{j\in J^+}\zeta_j(\rho) - e^{-i\rho}\sum_{j\in J^-}\zeta_j(\rho) = \DeltaP(\rho)\left[ \hat{q}_0(\rho) + \left( \frac{1}{\DeltaP(\rho)}-1 \right)H(\rho) \right],
\EE
where
\BES
H(\rho) = \sum_{j\in\widehat{J}^+}c_j(\rho)\tilde{h}_{\widehat{J}^+_j}(\rho) - e^{-i\rho}\sum_{j\in\widehat{J}^-}c_j(\rho)\tilde{h}_{\widehat{J}^-_j}(\rho),
\EES
\end{lem}
The other integrals in equation~\eqref{eqn:P1:Reps:Ser:q.Integrals.at.zeros} are around the boundaries of discs and circular sectors centred at each zero of $\DeltaP$. Over the next paragraphs we combine and simplify these integrals to the desired form.

Consider $\sigma\in D^+$ such that $\DeltaP(\sigma)=0$. Then the fourth integral on the right hand side of equation~\eqref{eqn:P1:Reps:Ser:q.Integrals.at.zeros} includes
\BES
\int_{C(\sigma,\varepsilon)}e^{i\rho x+i\rho^n(T-t)}\sum_{j\in J^+}\frac{\eta_j(\rho)}{\DeltaP(\rho)}\d\rho = \int_{C(\sigma,\varepsilon)}\frac{e^{i\rho x-i\rho^nt}}{\DeltaP(\rho)}\sum_{j\in J^+}{\zeta_j(\rho)}\d\rho,
\EES
the equality being justified by the following lemma, whose proof appears at the end of the section.
\begin{lem} \label{lem:P1:Reps:zeta.eta.relation}
Let $\Pi(n,A,a,q_0,h)$ be well-posed. Then for each $j\in\{1,2,\dots,2n\}$, the functions
\BE\label{eqn:P1:Reps:lem.zeta.eta.relation:entire.function}
\sum_{j\in\ J^+}\frac{\zeta_j(\rho) - e^{a\rho^nT}\eta_j(\sigma)}{\DeltaP(\rho)},\qquad \sum_{j\in\ J^-}\frac{\zeta_j(\rho) - e^{a\rho^nT}\eta_j(\sigma)}{\DeltaP(\rho)}
\EE
are entire.
\end{lem}

Consider $\sigma\in(\partial D)\cap\mathbb{C}^+$ such that $\DeltaP(\sigma)=0$. Define $\Gamma^D=\partial(B(\sigma,\varepsilon)\cap D)$ and $\Gamma^E=\partial(B(\sigma,\varepsilon)\cap E)$. Then the second and fourth integrals on the right hand side of equation~\eqref{eqn:P1:Reps:Ser:q.Integrals.at.zeros} include
\begin{gather*}
\int_{\Gamma^E}\frac{e^{i\rho x-i\rho^nt}}{\DeltaP(\rho)}\sum_{j\in J^+}{\zeta_j(\rho)}\d\rho \quad \mbox{and} \\
\int_{\Gamma^D}e^{i\rho x+i\rho^n(T-t)}\sum_{j\in J^+}\frac{\eta_j(\rho)}{\DeltaP(\rho)}\d\rho = \int_{\Gamma^D}\frac{e^{i\rho x-i\rho^nt}}{\DeltaP(\rho)}\sum_{j\in J^+}{\zeta_j(\rho)}\d\rho,
\end{gather*}
respectively, by Lemma~\ref{lem:P1:Reps:zeta.eta.relation}. The sum of the above expressions is
\BES
\int_{C(\sigma,\varepsilon)}\frac{e^{i\rho x-i\rho^nt}}{\DeltaP(\rho)}\sum_{j\in J^+}{\zeta_j(\rho)}\d\rho.
\EES

Consider $0\neq\sigma\in\mathbb{R}$ such that $\DeltaP(\sigma)=0$. Define $\Gamma^D=\partial(B(\sigma,\varepsilon)\cap D)$ and let $\Gamma^E=\partial(B(\sigma,\varepsilon)\cap E)$. Then the fourth and fifth integrals on the right hand side of equation~\eqref{eqn:P1:Reps:Ser:q.Integrals.at.zeros} include
\begin{gather*}
\int_{\Gamma^D}e^{i\rho x+i\rho^n(T-t)}\sum_{j\in J^+}\frac{\eta_j(\rho)}{\DeltaP(\rho)}\d\rho = \int_{\Gamma^D}\frac{e^{i\rho x-i\rho^nt}}{\DeltaP(\rho)}\sum_{j\in J^+}{\zeta_j(\rho)}\d\rho \quad \mbox{and} \\
\int_{\Gamma^E}\frac{e^{i\rho (x-1)-i\rho^nt}}{\DeltaP(\rho)}\sum_{j\in J^-}{\zeta_j(\rho)}\d\rho,
\end{gather*}
respectively, by analyticity and Lemma~\ref{lem:P1:Reps:zeta.eta.relation}. The sum of the above expressions is
\BES
\int_{C(\sigma,\varepsilon)}\frac{e^{i\rho x-i\rho^nt}}{\DeltaP(\rho)}\d\rho\sum_{j\in J^+}{\zeta_j(\sigma)} - \int_{\Gamma^E}\frac{e^{i\rho x-i\rho^nt}}{\DeltaP(\rho)}\left(\sum_{j\in J^+}{\zeta_j(\rho)} - e^{-i\rho}\sum_{j\in J^-}{\zeta_j(\rho)} \right)\d\rho.
\EES

Similar calculations may be performed for $\sigma\in E^-, D^-, (\partial D)\cap\mathbb{C}^-,\{0\}$. Define the index set $K^{\mathbb{R}}\subset\mathbb{N}$ by $k\in K^{\mathbb{R}}$ if and only if $\sigma_k\in\mathbb{R}$. For each $k\in K^\mathbb{R}$ define $\Gamma_k=\partial(B(\sigma_k,\varepsilon)\cap\mathbb{C}^-)$. Then, substituting the calculations above and applying Lemma~\ref{lem:P1:Reps:q0}, equation~\eqref{eqn:P1:Reps:Ser:q.Integrals.at.zeros} yields
\begin{multline*}
2\pi q(x,t) = \sum_{k\in K^+}\int_{C(\sigma_k,\varepsilon)}\frac{e^{i\rho x-i\rho^nt}}{\DeltaP(\rho)} \sum_{j\in J^+} \zeta_j(\rho)\d\rho \\ + \sum_{k\in K^-}\int_{C(\sigma_k,\varepsilon)}\frac{e^{i\rho(x-1)-i\rho^nt}}{\DeltaP(\rho)} \sum_{j\in J^-} \zeta_j(\rho)\d\rho \\
- \left\{\sum_{k\in K^\mathbb{R}} \int_{\Gamma_k^-} + \int_{\mathbb{R}} \right\} e^{i\rho x-i\rho^nt} \left( \frac{1}{\DeltaP(\rho)}-1 \right)H(\rho) \d\rho.
\end{multline*}
A residue calculation at each $\sigma_k$ completes the proof.

\subsection{Integral Representation} \label{ssec:P1:Reps:Int}
As $\Pi$ is well-posed, equation~\eqref{eqn:P1:Reps:Ser:q.Implicit} holds but, as $\Pi(n,A,-a,h,q_0)$ may not be well-posed, it is not possible to use Jordan's Lemma to deform the second and fifth integrals on the right hand side over $\widetilde{E}$. However it is still possible to deform the fourth and seventh integrals over $\widetilde{D}$. Hence two additional terms appear in equation~\eqref{eqn:P1:Reps:Ser:q.Integrals.at.zeros},
\BES
\int_{\partial\widetilde{E}^+}e^{i\rho x-i\rho^nt}\sum_{j\in J^+}\frac{\zeta_j(\rho)}{\DeltaP(\rho)}\d\rho + \int_{\partial\widetilde{E}^-}e^{i\rho(x-1)-i\rho^nt}\sum_{j\in J^-}\frac{\zeta_j(\rho)}{\DeltaP(\rho)}\d\rho.
\EES
The remainder of the derivation is unchanged from that presented in Section~\ref{ssec:P1:Reps:Ser}.

\subsection{Proofs of technical lemmata}

\begin{proof}[Proof of Lemma~\textup{\ref{lem:P1:Reps:q0}.}]
We expand the left hand side of equation~\eqref{eqn:P1:Repd:lem.q0:Statement} in terms of $u(\rho,l)$ and rearrange the result. To this end we define the matrix-valued function $X^{l\hspace{0.5mm}j}:\mathbb{C}\to\mathbb{C}^{(n-1)\times(n-1)}$ to be the $(n-1)\times(n-1)$ submatrix of
\BES
\BP\mathcal{A}&\mathcal{A}\\\mathcal{A}&\mathcal{A}\end{pmatrix}
\EES
whose $(1,1)$ entry is the $(l+1,r+j)$ entry. Then
\BE \label{eqn:P1:Repd:lem.q0:zetaj}
\widehat{\zeta}_j(\rho)=\sum_{l=1}^n u(\rho,l)\det X^{l\hspace{0.5mm}j}(\rho).
\EE

By Notation~\ref{ntn:P1:Implicit:DeltaP} and equation~\eqref{eqn:P1:Repd:lem.q0:zetaj}, the left hand side of equation~\eqref{eqn:P1:Repd:lem.q0:Statement} is equal to
\begin{multline*}
\sum_{l=1}^n u(\rho,l)\left[\left( \sum_{j:J_j\text{ odd}} c_{(J_j-1)/2}(\rho)\det X^{l\hspace{0.5mm}j} - \sum_{j:J'_j\text{ odd}} c_{(J'_j-1)/2}(\rho)\sum_{k=1}^n\M{\widehat{A}}{j}{k}\det X^{l\hspace{0.5mm}j} \right)\right. \\
\left.-e^{-i\rho}\left( \sum_{j:J_j\text{ even}} c_{J_j/2}(\rho)\det X^{l\hspace{0.5mm}j} - \sum_{j:J'_j\text{ even}} c_{J'_j/2}(\rho)\sum_{k=1}^n\M{\widehat{A}}{j}{k}\det X^{l\hspace{0.5mm}j} \right)\right] + H(\rho).
\end{multline*}
Splitting the sums over $k$ into $k:J_k$ is odd and $k:J_k$ is even and rearranging inside the parentheses, we evaluate the square bracket to
\BE \label{eqn:P1:Repd:lem.q0:Square.Bracket} \begin{split}
&\left[\sum_{j:J_j\text{ odd}} \left( c_{(J_j-1)/2}(\rho) - \sum_{\mathclap{k:J'_k\text{ odd}}} c_{(J'_k-1)/2}(\rho)\M{\widehat{A}}{k}{j} + e^{-i\rho}\sum_{\mathclap{k:J'_k\text{ even}}} c_{J'_k/2}(\rho)\M{\widehat{A}}{k}{j} \right)\det X^{l\hspace{0.5mm}j}\right. \\
&\hspace{3mm} \left.+ \sum_{j:J_j\text{ even}} \left( -c_{J_j/2}(\rho)e^{-i\rho} - \sum_{\mathclap{k:J'_k\text{ odd}}} c_{(J'_k-1)/2}(\rho)\M{\widehat{A}}{k}{j} + e^{-i\rho}\sum_{\mathclap{k:J'_k\text{ even}}} c_{J'_k/2}(\rho)\M{\widehat{A}}{k}{j} \right)\det X^{l\hspace{0.5mm}j}\right]\hspace{-1mm}.\end{split}
\EE
Making the change of variables $k\mapsto r$ defined by
\begin{gather*}
J'_k \mbox{ is odd if and only if } \widehat{J}^+\ni r = (J'_k-1)/2 \mbox{, in which case } k=\widehat{J}^+_r, \\
J'_k \mbox{ is even if and only if } \widehat{J}^-\ni r = J'_k/2 \mbox{, in which case } k=\widehat{J}^-_r,
\end{gather*}
it is clear that each of the parentheses in expression~\eqref{eqn:P1:Repd:lem.q0:Square.Bracket} evaluates to $\M{\mathcal{A}}{1}{j}$. Hence
\BES
\sum_{j\in J^+}\zeta_j(\rho) - e^{-i\rho}\sum_{j\in J^-}\zeta_j(\rho) = \sum_{l=1}^n u(\rho,l) \sum_{j=1}^n \M{\mathcal{A}}{1}{j}(\rho)\det X^{l\hspace{0.5mm}j}(\rho) + H(\rho). \qedhere
\EES
\end{proof}

\begin{proof}[Proof of Lemma~\textup{\ref{lem:P1:Reps:zeta.eta.relation}.}]
The $t$-transforms of the boundary functions are entire, as are the monomials $c_j$, hence the sum of products of a $t$-transform and monomials $c_j$ is also entire. By equation~\eqref{eqn:P1:Implicit:Zeta.Eta.f.g} this establishes that expressions~\eqref{eqn:P1:Reps:lem.zeta.eta.relation:entire.function} are entire functions of $\rho$.
\end{proof}

\section{Alternative characterisations} \label{sec:P1:Alt}

In this section we discuss sufficient conditions for well-posedness of initial-boundary value problems and present a proof of Theorem~\ref{thm:P1:Intro:Even.Can.Deform.If.Well-posed}. These topics are unified by the arguments and notation used.

\subsection{Sufficient conditions for well-posedness} \label{ssec:P1:Alt:Non-Robin}

Throughout Section~\ref{ssec:P1:Alt:Non-Robin} we assume the boundary conditions are non-Robin. This simplifies the PDE characteristic matrix greatly, leading to corresponding simplifications in the arguments presented below. Nevertheless, we identify suprising counterexamples to the qualitative hypothesis `highly coupled boundary conditions lead to well-posed problems whose solutions may be expressed using series.'

We give the condition whose effects are of interest.

\begin{cond} \label{cond:P1:Alt:Non-Robin:First.Condition}
For $A$, a boundary coefficient matrix specifying non-Robin boundary conditions, we define

$C=|\{j:\M{\alpha}{k}{j},\M{\beta}{k}{j}\neq0$ for some $k\}|$, the number of boundary conditions that couple the ends of the space interval, and

$R=|\{j:\M{\beta}{k}{j}=0$ for all $k\}|$, the number of right-handed boundary functions, whose corresponding column in $A$ is $0$.

Let $a=\pm i$ and let $A$ be such that
\BES
R \leq \left\{ \begin{matrix} \frac{n}{2} & \mbox{if } n \mbox{ is even and } a=\pm i \\ \frac{n+1}{2} & \mbox{if } n \mbox{ is odd and } a=i \\ \frac{n-1}{2} & \mbox{if } n \mbox{ is odd and } a=-i \end{matrix} \right\} \leq R+C.
\EES
\end{cond}

We investigate the effect of Condition~\ref{cond:P1:Alt:Non-Robin:First.Condition} upon the behaviour of the ratio
\BE \label{eqn:P1:Alt:Non-Robin:ratio}
\frac{\eta_m(\rho)}{\DeltaP(\rho)}
\EE
in the limit $\rho\to\infty$ from within $\widetilde{D}$. The PDE characteristic determinant is an exponential polynomial, a sum of terms of the form
\BES
Z(\rho)e^{-i\rho\sum_{y\in Y}\omega^{y}}
\EES
where $Z$ is some monomial and $Y\subset\{0,1,\dots,n-1\}$. As the problem may be ill-posed $\eta_m$ is defined as in Lemma~\ref{lem:P1:WellPosed:Ass.D.Implies.Admissible}, a sum of terms of the form
\BE \label{eqn:P1:Alt:Non-Robin:Term.of.etam}
X(\rho)e^{-i\rho\sum_{y\in Y}\omega^{y}}\int_0^1e^{-i\rho x\omega^{z}}q_T(x)\d x
\EE
where $X$ is some monomial, $q_T\in C^\infty[0,1]$, $Y\subset\{0,1,\dots,n-1\}$ and $z\in\{0,1,\dots,n-1\}\setminus Y$.

Fix $j\in\{1,2,\dots,n\}$ and let $\rho\in \widetilde{D}_j$. Then the modulus of
\BE \label{eqn:P1:Alt:Non-Robin:ExponentialY}
e^{-i\sum_{y\in Y}\omega^y\rho}
\EE
is uniquely maximised for the index set
\BES
Y=\begin{cases} \{j-1,j,\dots,j-2+\tfrac{n}{2}\} & n \mbox{ even,} \\ \{j-1,j,\dots,j-2+\tfrac{1}{2}(n+\Im(a))\} & n \mbox{ odd.} \end{cases}
\EES
By Condition~\ref{cond:P1:Alt:Non-Robin:First.Condition} $\DeltaP(\rho)$ has a term given by that exponential multiplied by some monomial coefficient, $Z_j(\rho)$. That term dominates all other terms in $\DeltaP(\rho)$ but it also dominates all terms in $\eta_m(\rho)$. Hence the ratio~\eqref{eqn:P1:Alt:Non-Robin:ratio} is bounded in $\widetilde{D}_j$ for each $j\in\{1,2,\dots,n\}$ and decaying as $\rho\to\infty$ from within $\widetilde{D}_j$.

If it were possible to guarantee that $Z_j\neq0$ then it would be proven that Condition~\ref{cond:P1:Alt:Non-Robin:First.Condition} is sufficient for well-posedness. Unfortunately this is not the case, as the following example shows.

\begin{eg} \label{eg:P1:Alt:Non-Robin:3Pseudo}
Let
\BES 
A=\BP1&-1&0&0&0&0\\0&0&1&-1&0&0\\0&0&0&0&1&2\EP
\EES
and consider the problem $\Pi(3,A,i,0,q_0)$. Then
\BES
\widetilde{D}_1\subseteq\left\{\rho\in\mathbb{C}:0<\arg\rho<\frac{\pi}{3}\right\}
\EES
and
\BES
\mathcal{A}(\rho) = \BP -c_2(\rho)(e^{-i\rho}-1) & -c_1(\rho)(e^{-i\rho}-1) & -c_0(\rho)(e^{-i\rho}+2) \\ -c_2(\rho)(e^{-i\omega\rho}-1) & -\omega c_1(\rho)(e^{-i\omega\rho}-1) & -\omega^2 c_0(\rho)(e^{-i\omega\rho}+2) \\ -c_2(\rho)(e^{-i\omega^2\rho}-1) & -\omega^2 c_1(\rho)(e^{-i\omega^2\rho}-1) & -\omega c_0(\rho)(e^{-i\omega^2\rho}+2) \EP.
\EES
We calculate
\begin{align*}
\DeltaP(\rho) &= (\omega-\omega^2)c_2(\rho)c_1(\rho)c_0(\rho)\left[9+(2-2)(e^{i\rho}+e^{i\omega\rho}+e^{i\omega^2\rho})\right. \\ &\hspace{60mm} \left.+ (1-4)(e^{-i\rho}+e^{-i\omega\rho}+e^{-i\omega^2\rho})\right], \\
\intertext{in this case, as $\M{\beta}{1}{2} + \M{\beta}{2}{1} + \M{\beta}{3}{0}=0$, the coefficients of $e^{i\omega^j\rho}$ cancel for each $j$, }
&= 3(\omega-\omega^2)c_2(\rho)c_1(\rho)c_0(\rho)\left[3-(e^{-i\rho}+e^{-i\omega\rho}+e^{-i\omega^2\rho})\right], \\
\eta_3(\rho)  &= (\omega^2-\omega)c_2(\rho)c_1(\rho)c_0(\rho)\sum_{j=0}^2{\omega^j\hat{q}_T(\omega^j\rho)(e^{i\omega^j\rho} - e^{-i\omega^{j+1}\rho} - e^{-i\omega^{j+2}\rho} + 1)}.
\end{align*}
Fix $\delta>0$. Consider a sequence, $(\rho_j)_{j\in\mathbb{N}}$, defined by $\rho_j=R_je^{i\pi/12}$ where
\BE \label{eqn:P1:Alt:Non-Robin:3Pseudo.eg:Bounds.On.Rj}
R_j\not\in\bigcup_{m=0}^\infty{\left(\left(1-\tfrac{\sqrt{3}}{3}\right)\sqrt{2}\pi m - \delta , \left(1-\tfrac{\sqrt{3}}{3}\right)\sqrt{2}\pi m + \delta\right)}
\EE
is a strictly increasing sequence of positive real numbers with limit $\infty$ chosen such that $\rho_j\in\widetilde{D}_1$ and. The ratio $\eta_3(\rho_j)/\DeltaP(\rho_j)$ evaluates to
\BES
\frac{ -\hat{q}_T(\rho_j) - \omega\hat{q}_T(\omega\rho_j)e^{-i(1-\omega)\rho_j} + \omega^2\hat{q}_T(\omega^2\rho_j)e^{i(\omega^2+\omega)\rho_j} + O(1)}{3(e^{-i(1-\omega)\rho_j} + 1) + O(e^{-R_j(\sqrt{3}-1)/2\sqrt{2}})}.
\EES
The denominator is $O(1)$ but, by condition~\eqref{eqn:P1:Alt:Non-Robin:3Pseudo.eg:Bounds.On.Rj}, is bounded away from $0$. The terms in the numerator all approach infinity at different rates, depending upon $\hat{q}_T$. Hence the ratio is unbounded and, by Theorem~\ref{thm:P1:Intro:WellPosed}, the problem is ill-posed.
\end{eg}

Indeed, third order initial-boundary value problems with pseudo-periodic boundary conditions are ill-posed if and only if
\begin{alignat*}{4}
a &=  i & &\mbox{ and } & \M{\beta}{1}{2} + \M{\beta}{2}{1} + \M{\beta}{3}{0} &= 0 & &\mbox{ or} \\
a &= -i & &\mbox{ and } & \frac{1}{\M{\beta}{1}{2}} + \frac{1}{\M{\beta}{2}{1}} + \frac{1}{\M{\beta}{3}{0}} &= 0. & &
\end{alignat*}
A combinatorial necessary and sufficient condition for $Z_j\neq0$ in odd order problems is presented as Condition~3.22 of~\cite{Smi2011a} but is omitted here due to its technicality however we do improve upon that condition; see Remark~\ref{rmk:P1:Spectrum:Symmetry.In.Coefficients}. No further third order examples are known which obey Condition~\ref{cond:P1:Alt:Non-Robin:First.Condition} but are ill-posed.

Condition~3.22 of~\cite{Smi2011a} may be adapted to even problems by setting $k=n/2-R$. The pseudo-periodic problems of second and fourth order are ill-posed if and only if
\begin{alignat*}{3}
n &= 2 & &\mbox{ and } & 0 &= \M{\beta}{1}{1} + \M{\beta}{2}{0}, \\
n &= 4 & &\mbox{ and } & 0 &= \M{\beta}{1}{3}\M{\beta}{2}{2} + \M{\beta}{2}{2}\M{\beta}{3}{1} + \M{\beta}{3}{1}\M{\beta}{4}{0} + \M{\beta}{4}{0}\M{\beta}{1}{3} + 2(\M{\beta}{1}{3}\M{\beta}{3}{1} + \M{\beta}{2}{2}\M{\beta}{4}{0}).
\end{alignat*}
For example, the problem $\Pi(4,A,\pm i,h,q_0)$ with boundary coefficient matrix
\BES
A=\BP1&1&0&0&0&0&0&0\\0&0&1&-1&0&0&0&0\\0&0&0&0&1&1&0&0\\0&0&0&0&0&0&1&-1\EP
\EES
is ill-posed.

\begin{rmk} \label{rmk:P1:Alt:Well-Posed.No.Series}
The essential difference between the odd and even cases presented above is that for odd order problems the well-posedness criteria depend upon the direction coefficient whereas for even order problems they do not. This means it is possible to construct examples of odd order problems that are well-posed but whose solutions cannot be represented by a series using Theorem~\ref{thm:P1:Intro:Reps.Ser}. Indeed the problem $\Pi(3,A,i,h,q_0)$, with boundary coefficient matrix given by
\BES
A=\BP1&-1&0&0&0&0\\0&0&1&-1&0&0\\0&0&0&0&1&\frac{1}{2}\EP,
\EES
is well-posed but is ill-posed in the opposite direction. This is the issue mentioned in Remark~3.3 of~\cite{FP2000a}.
\end{rmk}

\begin{rmk} \label{rmk:P1:Alt:Simple}
There are classes of examples for which $Z_j\neq0$ is guaranteed. Indeed, Condition~\ref{cond:P1:Alt:Non-Robin:First.Condition} is precisely the necessary and sufficient condition for well-posedness of problems with simple boundary conditions proved in~\cite{Pel2004a}.
\end{rmk}

\begin{rmk} \label{rmk:P1:Alt:Deformation.In.Upper.Half-Plane.Only}
There exist problems $\Pi$ for which $\Pi'$ is ill-posed but for which
\BES
\frac{\zeta_j(\rho)}{\DeltaP(\rho)}\to0 \mbox{ as } \rho\to\infty \mbox{ from within } \widetilde{E}^+,
\EES
for all $j\in J^+$ or from within $\widetilde{E}^-$ for all $j\in J^-$. This is a property of the $\zeta_j$, dependent upon which column of $\mathcal{A}$ is replaced with the transformed data, not of the sectors in which the decay or blow-up occurs. In this case it is possible to deform contours over the corresponding $\widetilde{E}^\pm$ hence one of the terms
\BES
\int_{\partial\widetilde{E}^+}e^{i\rho x-i\rho^nt}\sum_{j\in J^+}\frac{\zeta_j(\rho)}{\DeltaP(\rho)}\d\rho, \qquad \int_{\partial\widetilde{E}^-}e^{i\rho(x-1)-i\rho^nt}\sum_{j\in J^-}\frac{\zeta_j(\rho)}{\DeltaP(\rho)}\d\rho
\EES
evaluates to zero in equation~\eqref{eqn:P1:Intro:thm.Reps.Int:q} but the other does not.
\end{rmk}

\begin{rmk} \label{rmk:P1:Alt:Conjecture.Conditions.Necessary}
It is a conjecture that Condition~\ref{cond:P1:Alt:Non-Robin:First.Condition} together with Condition~3.22 of~\cite{Smi2011a} (as modified above to include $n$ even) are necessary as well as sufficient for well-posedness of problems with non-Robin boundary conditions. Any counterexample must satisfy several strong symmetry conditions that appear to be mutually exclusive. Indeed for a problem, which fails Condition~\ref{cond:P1:Alt:Non-Robin:First.Condition} or which satisfies Condition~\ref{cond:P1:Alt:Non-Robin:First.Condition} but for which $Z_j=0$, to be well-posed several monomial coefficients $X$ from equation~\eqref{eqn:P1:Alt:Non-Robin:Term.of.etam} must be identically zero.
\end{rmk}

\begin{rmk} \label{rmk:P1:Alt:Condition.Robin}
We give a condition equivalent to Condition~\ref{cond:P1:Alt:Non-Robin:First.Condition} for Robin type boundary conditions. Indeed, we define

$B_1=|\{j\in\widetilde{J}^-:\hsexists k,r$ for which $\M{\beta}{k}{j}\neq0$ and $\M{\alpha}{k}{r}$ is a pivot$\}|$,

$B_2=|\widetilde{J}^-|$ and

$B_3=|\{j\in\widetilde{J}^+:\hsexists k,r$ for which $\M{\alpha}{k}{j}\neq0$ and $\M{\beta}{k}{r}$ is a pivot$\}|$.

Then the condition is
\BES
B_2-B_1 \leq \left\{ \begin{matrix} \frac{n}{2} & \mbox{if } n \mbox{ is even and } a=\pm i \\ \frac{n+1}{2} & \mbox{if } n \mbox{ is odd and } a=i \\ \frac{n-1}{2} & \mbox{if } n \mbox{ is odd and } a=-i \end{matrix} \right\} \leq B_2+B_3.
\EES
\end{rmk}

\subsection{Series representations for $n$ even} \label{ssec:P1:Alt:Series.If.Even}

\begin{proof}[Proof of Theorem~\textup{\ref{thm:P1:Intro:Even.Can.Deform.If.Well-posed}.}]
By Theorem~\ref{thm:P1:Intro:WellPosed}, the well-posedness of $\Pi(n,A,i,h,q_0)$ and the arguments of Section~\ref{ssec:P1:Alt:Non-Robin}, for each $j\in\{1,2,\dots,n\}$ there exists some $\Ymax\subset\{0,1,\dots,n-1\}$ such that
\begin{enumerate}
\item{the term
\BES
Z_{\Ymax}(\rho)e^{-i\rho\sum_{y\in\Ymax}\omega^y}
\EES
appears in $\DeltaP$ with $Z_{\Ymax}\neq0$ a polynomial,
\BES
\DeltaP(\rho)=O\left(|Z_{\Ymax}(\rho)|e^{\Im\left(\rho\sum_{y\in\Ymax}\omega^y\right)}\right) \mbox{ as } \rho\to\infty \mbox{ from within } \widetilde{D}_j
\EES}
\item{for all $Y\subset\{0,1,\dots,n-1\}$, $z\in\{0,1,\dots,n-1\}\setminus Y$ for which
\BES
\M{X}{Y}{z}(\rho)e^{-i\rho\sum{y\in Y}\omega^y}\hat{q}_T(\omega^z\rho)
\EES
is a term in $\eta_k$, for some $k$, with $\M{X}{Y}{z}\neq0$ a polynomial such that
\BES
\M{X}{Y}{z}(\rho)e^{-i\rho\sum{y\in Y}\omega^y}\hat{q}_T(\omega^z\rho) = o\left(|Z_{\Ymax}(\rho)|e^{\Im\left(\rho\sum_{y\in\Ymax}\omega^y\right)}\right)
\EES
as $\rho\to\infty$ from within $\widetilde{D}_j$.}
\end{enumerate}
Hence, for all such $Y$, $z$,
\BE \label{eqn:P1:Alt:Series.If.Even:Well-posed.Condition}
\Im\left[e^{i\phi}\left(\sum_{y\in Y}\omega^y+\omega^z-\sum_{y\in\Ymax}\omega^y\right)\right] < 0
\EE
for all $x\in(0,1)$ and for all $\phi\in(0,\pi/n)$.

If $\Pi(n,A,-i,h,q_0)$ is ill-posed then there exist $Y$, $z$ satisfying the conditions above, $x\in(0,1)$ and $\phi\in(\pi/n,2\pi/n)$ such that
\BES
N = \Im\left[e^{i\phi}\left(\sum_{y\in Y}\omega^y+\omega^z-\sum_{y\in\Ymax}\omega^y\right)\right] > 0
\EES
Define $\overline{Y}_{\mathrm{max}}=\{n-y:y\in\Ymax\}$ and
\BES
\overline{Y} = \begin{cases}Y\cup\{z\} & \Im(e^{i\phi}\omega^z)\geq0, \\ Y & \Im(e^{i\phi}\omega^z)<0.\end{cases}
\EES
Then, as $n$ is even,
\BES
\Im\left[e^{i\phi}\left(\sum_{y\in \overline{Y}}\omega^y+\sum_{y\in\overline{Y}_{\mathrm{max}}}\omega^y\right)\right] \geq N > 0,
\EES
hence there exists some $\bar{x}\in(0,1)$ such that
\BES
\Im\left[e^{i(\phi-\frac{\pi}{n})}\left(\sum_{y\in Y}\omega^y+\bar{x}\omega^z-\sum_{y\in\Ymax}\omega^y\right)\right] > 0,
\EES
which contradicts inequality~\eqref{eqn:P1:Alt:Series.If.Even:Well-posed.Condition}. The argument is identical in the other direction, switching the intervals in which $\phi$ lies.
\end{proof}

\section{PDE discrete spectrum} \label{sec:P1:Spectrum}

In this section we investigate the PDE discrete spectrum, the set of zeros of an exponential polynomial. We use the definitions, results and arguments presented in~\cite{Lan1931a}.

\begin{lem} \label{lem:P1:Spectrum:Properties}
The PDE characteristic determinant and PDE discrete spectrum have the following properties:
\begin{enumerate}
\item{$\DeltaP(\rho)=(-1)^{n-1}\DeltaP(\omega\rho)$.}
\item{let $Y\subset\{0,1,\dots,n-1\}$, $Y'=\{y+1\mod n:y\in Y\}$. Let $Z_Y$ and $Z_{Y'}$ be the polynomial coefficients of $\exp{(-i\rho\sum_{y\in Y}\omega^y)}$ and $\exp{(-i\rho\sum_{y\in Y'}\omega^y)}$, respectively, in $\DeltaP$. Then $Z_Y(\rho) = (-1)^{n-1}Z_{Y'}(\omega\rho)$.}
\item{either $\DeltaP$ is a polynomial or the PDE discrete spectrum is asymptotically distributed in finite-width semi-strips each parallel to the outward normal to a side of a polygon with order of rotational symmetry a multiple of $n$. Further, the radial distribution of the zeros within each strip is asymptotically inversely proportional to the length of the corresponding side.}
\end{enumerate}
\end{lem}

\begin{proof}
(i) The identity
\BES
\M{\mathcal{A}}{k}{j}(\omega\rho) = \M{\mathcal{A}}{k+1}{j}(\rho)
\EES
follows directly from the definition~\eqref{eqn:P1:Intro:PDE.Characteristic.Matrix} of the PDE characteristic matrix. A composition with the cyclic permutation of order $n$ in the definition of the determinant yields the result.

(ii) By definition there exist a collection of index sets $\mathcal{Y}\subset\mathcal{P}\{0,1,\dots,n-1\}$ and polynomial coefficients $Z_Y(\rho)$ such that
\BES
\DeltaP(\rho) = \sum_{Y\in\mathcal{Y}}Z_Y(\rho)e^{-i\rho\sum_{y\in Y}\omega^y}.
\EES
By part (i),
\BES
\sum_{Y\in\mathcal{Y}}Z_Y(\rho)e^{-i\rho\sum_{y\in Y}\omega^y} = (-1)^{n-1}\sum_{Y\in\mathcal{Y}}Z_Y(\omega\rho)e^{-i\rho\sum_{y\in Y}\omega^{y+1}}
\EES
Define the collection $\mathcal{Y}'=\{\{y+1\mod n:y\in Y\}:Y\in\mathcal{Y}\}$. Then
\BES
\sum_{Y\in\mathcal{Y}}Z_Y(\rho)e^{-i\rho\sum_{y\in Y}\omega^y} = (-1)^{n-1}\sum_{Y'\in\mathcal{Y}'}Z_{Y'}(\rho)e^{-i\rho\sum_{y\in Y'}\omega^y}.
\EES
Equating coefficients of $\exp{(-i\rho\sum_{y\in Y}\omega^y)}$ yields $\mathcal{Y}=\mathcal{Y}'$ and the result follows.

(iii) The result follows from part (ii) and Theorem~8 of~\cite{Lan1931a}.
\end{proof}

An immediate corollary of Lemma~\ref{lem:P1:Spectrum:Properties} is that the PDE discrete spectrum has no finite accumulation point and is separated by some $\varepsilon>0$.

\begin{rmk} \label{rmk:P1:Spectrum:Symmetry.In.Coefficients}
A corollary of (ii) is that $Z_Y=0$ if and only if $Z_{Y'}=0$. This means it is only necessary to check $Z_j\neq0$ for a particular $j$ in conjunction with Condition~\ref{cond:P1:Alt:Non-Robin:First.Condition} to ensure well-posedness. This permits a simplification of the general Condition~3.22 of~\cite{Smi2011a}.
\end{rmk}

It is possible to strengthen part (iii) of Lemma~\ref{lem:P1:Spectrum:Properties} in certain cases.

\begin{thm} \label{thm:P1:Spectrum:Odd.Rays}
Let $n\geq3$ be odd and let $A$ be such that $\DeltaP$ is not a polynomial. If $n\geq 7$ we additionally require that Condition~\ref{cond:P1:Alt:Non-Robin:First.Condition} holds and the relevant coefficients, $Z_j$, are all nonzero. Then the PDE discrete spectrum must lie asymptotically on rays instead of semi-strips.
\end{thm}

\begin{proof}
Assume $n\geq7$ and the additional conditions hold. If
\BES
Y=\left\{1,2,\dots,\frac{n-1}{2}\right\}\in\mathcal{Y}
\EES
then, by part (ii) of Lemma~\ref{lem:P1:Spectrum:Properties}, $\{0,1,\dots,(n-3)/2\}\in\mathcal{Y}$ hence the indicator diagram of $\DeltaP$ has subset the convex hull of
\BES
S=\left\{\overline{\omega^r\sum_{y\in Y}\omega^y}:r\in\{0,1,\dots,n-1\}\right\}.
\EES
If $\overline{Y}=\{1,2,\dots,(n+1)/2\}\in\mathcal{Y}$ then the indicator diagram contains the regular $2n$-gon that forms the convex hull of $S\cup\{-s:s\in S\}$. We show that the indicator diagram is precisely the convex hull of $S$ or of $S\cup\{-s:s\in S\}$ and that there are no points $\sum_{y\in Y'}\omega^y$, for $Y'\in\mathcal{Y}$, on the boundary of the indicator diagram other than at the vertices.

Excepting rotations of $Y$ and $\overline{Y}$, which all correspond to vertices, the sets $Y'\in\mathcal{Y}$ whose corresponding exponent has greatest modulus
\BES
s'=\left|\sum_{y\in Y'}\omega^y\right|
\EES
are rotations and reflections of
\begin{align*}
Y_1 &= \{1,2,\dots,(n-3)/2\} \mbox{ or } \\
Y_2 &= \{1,2,\dots,(n-3)/2,(n+1)/2\}.
\end{align*}
However, the minimum modulus of the boundary of the indicator diagram is greater than or equal to
\begin{multline*}
\frac{1}{2} \left| \sum_{y=0}^{(n-3)/2}\omega^y + \sum_{y=1}^{(n-1)/2}\omega^y \right| = \left| \frac{1}{2}(1+\omega^{(n-1)/2}) + \sum_{y=1}^{(n-3)/2}\omega^y \right| \\
> \left| \sum_{y=1}^{(n-3)/2}\omega^y \right| = s_1 > \left| \omega^{(n+1)/2} + \sum_{y=1}^{(n-3)/2}\omega^y \right| = s_2,
\end{multline*}
hence any point corresponding to $Y'$ is interior to the indicator diagram. It is easy to check that this also holds if $n=3$ or $n=5$.

As there can only be two colinear exponents lying on any side of the indicator diagram, the argument in Sections~1--7 of~\cite{Lan1931a} may be simplified considerably and yield the stronger condition that the zeros of the exponential polynomial lie asymptotically on a ray, a semi-strip of zero width. The arguments of Sections~8--9 applied to this result complete the proof.
\end{proof}

\begin{rmk} \label{rmk:P1:Spectrum:Rays.2.4.7}
Theorem~\ref{thm:P1:Spectrum:Odd.Rays} does not hold for $n$ even. Indeed,
\BES
\frac{1}{2}\left(\sum_{j=0}^{(n-2)/2}\omega^j - \sum_{j=1}^{n/2}\omega^j  \right) = \sum_{j=1}^{(n-2)/2}\omega^j,
\EES
hence if
\BE \label{eqn:P1:Spectrum:rmk.Rays.2.4.7:Even}
\{0,1,\dots,(n-2)/2\},\{1,2,\dots,n/2\},\{1,2,\dots,(n-2)/2\} \in \mathcal{Y}
\EE
then, by part 2.\ of Lemma~\ref{lem:P1:Spectrum:Properties}, there are three colinear exponents on each side of the indicator diagram. Condition~\eqref{eqn:P1:Spectrum:rmk.Rays.2.4.7:Even} does not represent a pathological counterexample; it is satisfied by most pseudoperiodic, including all quasiperiodic, boundary conditions.
\end{rmk}

\bigskip

The author is sincerely grateful to B. Pelloni for her continued support and encouragement. He is funded by EPSRC.

\bibliographystyle{amsplain}
\bibliography{dbrefs}
\end{document}